\documentclass[11pt]{amsart}
\pdfoutput=1
\usepackage[utf8]{inputenc}
\usepackage[english]{babel}
\usepackage{attachfile}
\usepackage{navigator}
\usepackage{embedfile}

\embedfile{Coefficients_h_k=1.txt}
\embedfile{Coefficients_h_k=2.txt}
\embedfile{Coefficients_h_k=3.txt}
\embedfile{Coefficients_h_k=4.txt}
\embedfile{Coefficients_h_k=5.txt}
\embedfile{Coefficients_c_k=1.txt}
\embedfile{Coefficients_c_k=2.txt}
\embedfile{Coefficients_c_k=3.txt}
\embedfile{Coefficients_c_k=4.txt}
\embedfile{Coefficients_c_k=5.txt}

\usepackage{enumitem}
\usepackage[T1]{fontenc}
\usepackage{amssymb,amsfonts,stmaryrd,mathrsfs}
\usepackage{amsmath,latexsym,amsthm,mathtools,bbm}
\usepackage{tikz-cd}
\usepackage{relsize,exscale}
\usepackage{color}

 \usepackage[capitalize]{cleveref}
 \usepackage{graphicx,xspace}

\usepackage{bm,changes}
\usepackage{graphicx}
\usepackage[export]{adjustbox}
\usepackage{pgf,tikz} 
\usepackage{mathrsfs}
\usetikzlibrary{arrows}
\usepackage{multicol}
\usepackage{caption}
\usepackage{subcaption}

\theoremstyle{plain}
\newtheorem{thm}{Theorem}[section]
\newtheorem{cor}[thm]{Corollary}
\newtheorem{lem}[thm]{Lemma}
\newtheorem{prop}[thm]{Proposition}

\theoremstyle{definition}
\newtheorem{conj}{Conjecture}
\newtheorem{defi}{Definition}[section]
\newtheorem{defi-lem}[defi]{Definition-Lemmma}

\theoremstyle{remark}
\newtheorem{rmq}{Remark}
\newtheorem{exe}{Example}[section]

\DeclareMathOperator{\hook}{hook}
\DeclareMathOperator{\w}{w}

\definecolor{vert}{rgb}{0.2, 0.68, 0}

\author{Houcine Ben Dali}
\address{\parbox{\linewidth}{Université de Paris, CNRS, IRIF, F-75006, Paris, France \\
Université de Lorraine, CNRS, IECL, F-54000 Nancy, France}}

\email{bendali@irif.fr}

\begin{document}
\title[Generating series of non-oriented constellations]{Generating series of non-oriented constellations and marginal sums in the Matching-Jack conjecture}

\begin{abstract}
Using the description of hypermaps with matchings, Goulden and Jackson have given an expression of the generating series of rooted bipartite maps in terms of the zonal polynomials. We generalize this approach to the case of constellations on non-oriented surfaces that have recently been introduced  by Chapuy and Do\l{}{\k{e}}ga. A key step in the proof is an encoding of constellations with tuples of matchings.

We consider a one parameter deformation of the generating series of constellations using Jack polynomials and we introduce the coefficients $c^\lambda_{\mu^0,...,\mu^k}(b)$ obtained by the expansion of these functions in the power-sum basis. These coefficients are indexed by $k+2$ integer partitions and the deformation parameter $b$, and can be considered as a generalization for $k\geq1$ of the connection coefficients introduced by Goulden and Jackson. We prove that when we take some marginal sums, these coefficients enumerate $b$-weighted $k$-tuples of matchings. This can be seen as an "unconnected" version of a recent result of Chapuy and Do\l{}{\k{e}}ga for constellations. For $k=1$, this gives a partial answer to Goulden and Jackson Matching-Jack conjecture.

Lassale has formulated a positivity conjecture for the coefficients $\theta^{(\alpha)}_\mu(\lambda)$, defined as the coefficient of the Jack polynomial $J_\lambda^{(\alpha)}$ in the power-sum basis. We use the second main result of this paper to give a proof of this conjecture in the case of partitions $\lambda$ with rectangular shape.
\end{abstract}
\maketitle  

\section{Introduction}
\subsection{Jack polynomials and maps.}
Jack polynomials $J_\theta^{(\alpha)}$ are symmetric functions indexed by an integer partition $\theta$ and a deformation parameter $\alpha$ that were introduced in \cite{J}. Jack polynomials can be considered as one parameter deformation of Schur functions, which are obtained by evaluating the Jack polynomials at $\alpha=1$. For $\alpha=2$ we recover the zonal polynomials. This family of symmetric functions is related to various combinatorial problems \cite{Han88, GJ96b, DFS14}.
Some properties of Jack polynomials have been investigated in \cite{stan89} and \cite[Chapter VI]{Mac95}.

 In this paper, we will be interested in relationships between Jack polynomial series and generating series of maps. A connected map is a 2-cell embedding of a  connected graph into a closed surface without boundary, orientable or not. A map\footnote{This is not the usual definition of maps; what is usually called a map is called here a connected map.} is an unordered collection of connected maps. In this paper, we will use the word \textit{orientable} for maps on orientable surfaces and the word \textit{non-oriented} for maps on general surfaces, orientable or not. Maps appear in various branches of algebraic combinatorics, probability and physics. The study of maps involves various methods such as generating series, matrix integral techniques and bijective methods, see e.g  \cite{LZ04,Eyn16,BC86,Cha11,AL20}.
In this paper we will consider a class of vertex-colored maps that generalize bipartite maps, called $k$-constellations. Constellations on orientable surfaces were introduced in \cite{LZ04} and were generalized to the case of non-orientable surfaces in \cite{CD20}, see \Cref{ssMatchings}.

Let $j_\theta^{(\alpha)}:=\langle J_\theta,J_\theta\rangle_\alpha$ be the squared norm of the Jack polynomial associated to $\theta$ with respect to the $\alpha$-deformation of the Hall scalar product, see \cref{sec SymFun}. We consider $k+2$ different alphabets $\mathbf{x}^{(i)}:=(x^{(i)}_1,x^{(i)}_2,...)$, for $-1\leq i\leq k$, and we set the power-sum variables associated respectively to these alphabets $\mathbf{p}:=(p_1,p_2...)$ and $\mathbf{q}^{(i)}:=(q^{(i)}_1,q^{(i)}_2...)$ for $0\leq i\leq k$.
Chapuy and Do\l{}{\k{e}}ga have introduced\footnote{The function introduced in \cite{CD20} is a specialization of this function.}  in \cite{CD20} for every $k\geq1$ a function $\tau^{(k)}_b$ with $k+2$ sets of variables, defined as follows: 
\begin{equation}\label{eqtau}
    \tau^{(k)}_b(t,\mathbf{p},\mathbf{q}^{(0)},..,\mathbf{q}^{(k)}):=\sum_{n\geq0}t^n\sum_{\theta\vdash n}\frac{1}{j^{(\alpha)}_\theta}J^{(\alpha)}_\theta(\mathbf{p})J^{(\alpha)}_\theta(\mathbf{q}^{(0)})...J^{(\alpha)}_\theta(\mathbf{q}^{(k)}),
\end{equation}
where $J_\theta^{(\alpha)}$ are the Jack polynomials of parameter $\alpha=b+1$, expressed in the power-sum variables  $J_\theta^{(\alpha)}(\mathbf{p}):=J_\theta^{(\alpha)}(\mathbf{x}^{(-1)})$ and $J_\theta^{(\alpha)}(\mathbf{q}^{(i)}):=J_\theta^{(\alpha)}(\mathbf{x}^{(i)})$ for $0\leq i\leq k$. This function can be seen as an extension to $k+2$ sets of variables of the Cauchy sum for Jack symmetric functions.  We also define
\begin{equation}\label{eqPsi}
    \Psi^{(k)}_b(t,\mathbf{p},\mathbf{q}^{(0)},\mathbf{q}^{(1)},...,\mathbf{q}^{(k)}):=(1+b)t\frac{\partial}{\partial t}\log \tau^{(k)}_b(t,\mathbf{p},\mathbf{q}^{(0)},...,\mathbf{q}^{(k)}).
\end{equation}
In the case $k=1$, these functions were first introduced by Goulden and Jackson \cite{GJ96b}. They suggested that the function $\tau^{(1)}_b$ is related to generating series of matchings and $\Psi^{(1)}_b$ is related to the generating series of connected hypermaps (or by duality connected bipartite maps). The exponent of the shifted parameter $b:=\alpha-1$ is claimed to be correlated to the  bipartiness of the matchings in the first case and to the orientability of the maps in the second one. This was formulated in two conjectures that are still open, namely the $b$-conjecture and the Matching-Jack conjecture. These conjectures imply that the coefficients of the functions $\tau_b^{(1)}$ and $\Psi_b^{(1)}$ in the power-sums basis denoted respectively $c^\lambda_{\mu,\nu}(b)$ and $h^\lambda_{\mu,\nu}(b)$  are non-negative integer polynomials in $b$.
In  this paper, we consider a generalization of these quantities $c^\lambda_{\mu^0,...,\mu^k}$ and $h^\lambda_{\mu^0,..,\mu^k}$ indexed by $k+2$ partitions and defined by the expansion of $\tau^{(k)}_b$ and $\Psi_b^{(k)}$ in power-sum basis (see \cref{defc} and \cref{defh}).  We investigate their relationship with the enumeration of non-oriented constellations and tuples of matchings.

\subsection{Main Contributions.}
We now say a word about the main contributions of the paper. The first four points will be discussed in more details in the next subsections of the introduction:
\begin{itemize}
\item  We describe an encoding of non-oriented constellations with tuples of matchings; two versions of this correspondence are given, see \cref{prop1} and \cref{prop2}.
\item This correspondence is used to obtain \cref{Thm b=1} which relates the generating series of non-oriented $k$-constellations to the function $\Psi_b^{(k)}$ in the special case $b=1$. The case $k=1$ of this result was proved by Goulden and Jackson in \cite{GJ96a}.

\item In the second part of this paper, we consider some marginal sums of the coefficients $c^\lambda_{\mu^0,...,\mu^k}$ and $h^\lambda_{\mu^0,..,\mu^k}$, where we control two partitions $\lambda$ and $\mu$ and the number of parts of the other partitions, denoted respectively  $c^\lambda_{\mu,l_1,...l_k}(b)$ and $h^\lambda_{\mu, l_1,...,l_k}(b)$. \cref{Thm3} (see also \cref{thm marginal sums})  states that the coefficients $c^\lambda_{\mu,l_1,...l_k}(b)$ are non-negative integer polynomials in $b$ and that they enumerate $b$-weighted $k$-tuples of matchings. 
The proof is based on the work of Chapuy and Do\l{}{\k{e}}ga \cite{CD20} that gives an analog result for the coefficients $h^\lambda_{\mu, l_1,...,l_k}(b)$. The fact that the coefficients $c^\lambda_{\mu,l_1,...l_k}(b)$ are polynomials in $b$ with positive coefficients can directly be obtained from the result of \cite{CD20}, but not the integrality because of the derivative taken in \cref{eqPsi}. In the proof of \cref{Thm3} we use symmetry properties to eliminate factors appearing in the denominator. 
When $k=1$, \cref{Thm3} gives the marginal sum case in the Matching-Jack conjecture, and covers other partial results established for this conjecture \cite{KV16,KPV18}.

\item \Cref{Thm Lassale conjecture} gives a combinatorial expression for the coefficients of the development of Jack polynomials $J_\lambda$ in the power-sum basis, for rectangular partitions $\lambda$. In particular, this completes the proof of a Lassale's conjecture in the rectangular case. The proof is based on \Cref{Thm3}.

\item \cref{top degree 1} and \cref{top degree 2} give a combinatorial interpretation of the top degree part in coefficients $c^\lambda_{\mu^0,...,\mu^k}$. In the case $k=1$, this was investigated  in \cite{B21} using Jack characters, we give here a different proof.
\end{itemize}

\subsection{Constellations and matchings.}\label{ssMatchings}

We consider the definition of constellations on general surfaces, orientable or not, given in \cite{CD20}. The link with the usual definition of constellations in the orientable case is explained in \Cref{sec cons}.
\begin{defi}\label{def const}
Let $k\geq 1$. A (non-oriented) $k$-constellation is a map, connected or not, whose vertices are colored with colors $\{0,1,...,k\}$, such that\footnote{We use here the convention of \cite{CD20}, what we call $k$-constellation is often called $k+1$-constellation in the orientable case.}:
\begin{enumerate}
    \item Each vertex of color 0 (respectively $k$) has only neighbors of colors 1 (respectively $k-1$).
    \item For $0<i<k$, a vertex of color $i$ has only neighbors of color $i-1$ and $i+1$, and each corner of such vertex separates two vertices of colors $i-1$ and $i+1$.
\end{enumerate}
\end{defi}

\begin{figure}[t]
    \centering
    \includegraphics[width=.3\textwidth,center]{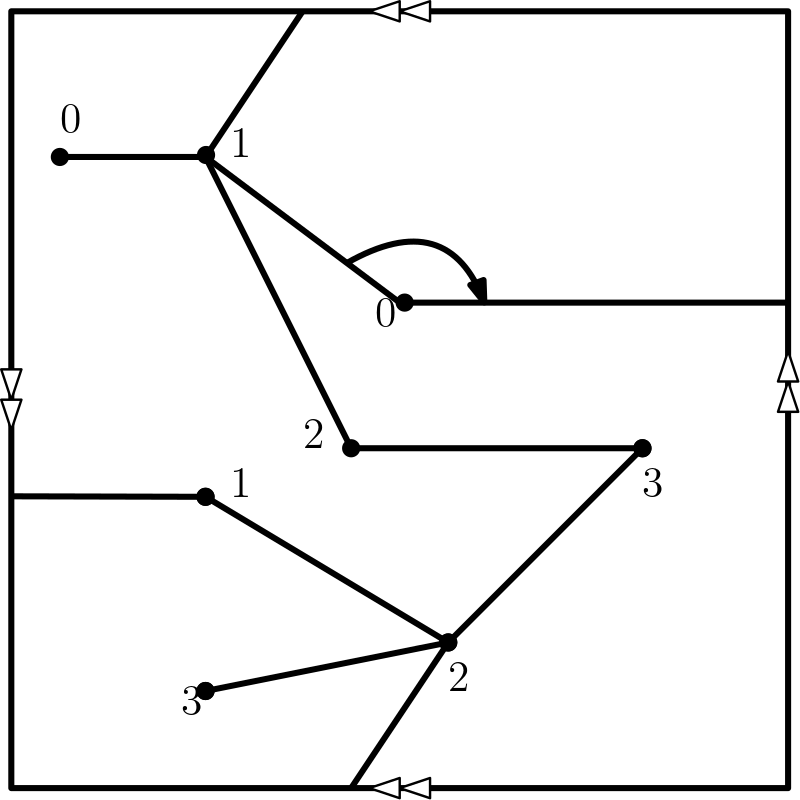}
    \caption{An example of a rooted 3-constellation drawn on the Klein Bottle. The
left-hand side of the square should be glued to the right-hand
one (with a twist) and the top side should be glued to the
bottom one (without a twist), as indicated by the white arrows. The root corner is indicated by a black arrow.
}
    \label{ex 3-cons Klein}
\end{figure}

Constellations come with a natural notion of rooting; a connected constellation $\mathbf{M}$ is \textit{rooted} by distinguishing an oriented corner $c$ of color 0. The rooted constellation obtained is denoted $(\mathbf{M},c)$. We can define the \textit{size} of a constellation as the sum of the degrees of all vertices of color 0. An example of a rooted non-oriented 3-constellation of size 3 is illustrated in \cref{ex 3-cons Klein}.
Given a $k$-constellation $\mathbf{M}$, we can define its \textit{profile} as the $k+2$-tuple of integer partitions $(\lambda,\mu^0,...,\mu^k)$, where $\lambda$ is the distribution of the face degrees, and $\mu^i$ is the distribution of the vertices of color $i$ for $0\leq i\leq k$, see \cref{sec cons} for a precise definition.

Constellations in the orientable case have been studied for a long time and have a well-known  combinatorial description given by rotation systems; an orientable $k$-constellation of size $n$ can be encoded by $k+1$-tuple of permutations of $\mathfrak{S}_n$. Moreover the profile of the constellation is related to the cyclic type of the permutations, see \cite{LZ04,BMS00,Cha18,F16}. In \cref{sec cons}, we introduce a notion of labelling for non-oriented constellations using \textit{right-paths}. This leads to a correspondence between $k$-constellations and $k+2$-tuples of matchings, see \cref{prop1}. In fact, this completes the following table;

\begin{center}
   \begin{tabular}{ c | c | c  }
     
      & bipartite maps & $k$-constellations \\ \hline
     orientable & pairs of permutations [classical]  & $k+1$-tuples of permutations \cite{JV90} \\\hline
     non-oriented & triples of matchings \cite{GJ96a} & $k+2$-tuples of matchings, this paper.
    
   \end{tabular}
 \end{center}

For $n\geq1$, we will consider matchings on the set $\mathcal{A}_n:=\{1,\hat{1},...,n,\hat{n}\}$, that is a set partition of $\mathcal{A}_n$ into pairs. A matching $\delta$ on $\mathcal{A}_n$ is \textit{bipartite} if each one of its pairs is of the form $(i,\hat j)$. We denote by $\varepsilon$ the bipartite matching on $\mathcal{A}_n$ defined by $\varepsilon:=\left\{\{1,\hat{1}\},\{2,\hat{2}\},...,\{n,\hat{n}\}\right\}$. Given two matchings $\delta_1$ and $\delta_2$, we consider the partition $\Lambda(\delta_1,\delta_2)$ of $n$ obtained by reordering the half-sizes of the connected components of the graph formed by $\delta_1$ and $\delta_2$.
For each $\lambda$ partition of $n$, we set the matching
\begin{align}\label{def delta lambda}
    \delta_\lambda:=&\big\{ \{1,\hat{2} \},\{2,\hat{3}\}.. \{\lambda_1-1,\widehat{\lambda_1} \}, \{\lambda_1,\hat{1}\},\\
    &\{\lambda_1+1,\widehat{\lambda_1+2}\},...\big\}\nonumber.
\end{align}

\noindent The matchings $\varepsilon$ and $\delta_\lambda$ will have a specific role in this article as an example of bipartite matchings satisfying $\Lambda(\varepsilon,\delta_\lambda)=\lambda$.

\begin{figure}[t]
    \centering
    \includegraphics[width=.4\textwidth, center]{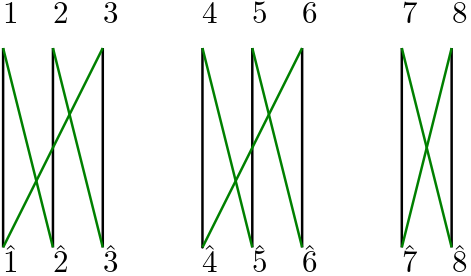}
    \caption{$G(\varepsilon,\delta_\lambda)$ for $\lambda=[3,3,2]$, $\varepsilon$ in black and $\delta_\lambda$ in green.}
    \label{fig matchings}
\end{figure}

\begin{exe}
In \cref{fig matchings}, the graph formed by $\varepsilon$ and $\delta_\lambda$ is illustrated for $n=8$ and $\lambda=[3,3,2]$.
\end{exe} Note that the symmetric group $\mathfrak{S}_n$ is in natural bijection with the set of bipartite matchings of $\mathcal{A}_n$  via the map $\sigma\longmapsto \delta_\sigma$ where $\delta_\sigma$ is the bipartite matching whose pairs are $(i,\widehat{\sigma(i)})_{1\leq i\leq n}.$ Given a $k$-tuple of matchings, $\lambda,\mu^0,...,\mu^k\vdash n$, we consider the two sets
\begin{multline}\label{mathfrakF}
 \mathfrak{F}^\lambda_{\mu^0,...,\mu^k}:=
     \Big\{(\delta_0,...,\delta_{k-1}) \text{ matchings of $\mathcal{A}_n$ such that } \Lambda(\varepsilon,\delta_0)=\mu^0\text{, }\Lambda(\delta_{k-1},\delta_{\lambda})=\mu^k \\
    \text{ and } \Lambda(\delta_{i-1},\delta_i)=\mu^i \hspace{0.3cm} \forall i \in\{1,...,k-1\}\Big\},
\end{multline}
and

 $$\tilde{\mathfrak{F}}^\lambda_{\mu^0,...,\mu^k}:=
     \Big\{(\delta_0,...,\delta_{k-1})\in\mathfrak{F}^\lambda_{\mu^0,...,\mu^k} \text{such that } \delta_i
    \text{ is bipartite  for } 1\leq i\leq k-1 \Big\}.$$

In \cref{prop2}, we prove that $k$-constellations with profile $(\lambda,\mu^0,...,\mu^k)$ with a specific labelling are in bijection with $ \mathfrak{F}^\lambda_{\mu^0,...,\mu^k}$. 

\subsection{Generating series of constellations.}\label{sec gs of const}
As described in \cref{ssMatchings}, constellations on the orientable case can be encoded with tuples of permutations. The multivariate generating series of permutations with respect to their cyclic type,  can be related via representation theory tools to the function $\tau^{(k)}_0$, see e.g \cite[Proposition 1.1]{GJ08} (in the following, we formulate this result using bipartite matchings rather than permutations). This gives a formula for the generating series of orientable constellations. To state this result we need to introduce the following notation; for a partition $\lambda$ and a non-negative integer $n$ we write $\lambda\vdash n$ if $\lambda$ is a partition of $n$ and if $\lambda=[\lambda_1,\lambda_2,..]$ we set as in \cite{Mac95}
$$z_\lambda:=\prod_{i\geq1}m_i(\lambda)!i^{m_i(\lambda)},$$ where $m_i(\lambda)$ is the number of parts of $\lambda$ equal to $i$. For a set of variables $\mathbf{u}:=(u_1,u_2,..)$, we denote $u_\lambda:=u_{\lambda_1}u_{\lambda_2}$...

\begin{thm}[\cite{GJ08}]\label{Thm b=0}
For $k\geq1$, we have
$$(i)\hspace{1cm} 
   \tau^{(k)}_0(t,\mathbf{p},\mathbf{q}^{(0)},..,\mathbf{q}^{(k)})=\sum_{n\geq 0}t^n\sum_{\lambda,\mu^0,...,\mu^k\vdash n}\frac{|\tilde{\mathfrak{F}}^\lambda_{\mu^0,...,\mu^k}|}{z_\lambda}p_\lambda q^{(0)}_{\mu^{0}}..q^{(k)}_{\mu^{k}},
$$
$$(ii)\hspace{1cm}\Psi_0^{(k)}(t,\mathbf{p},\mathbf{q}^{(0)},..,\mathbf{q}^{(k)})=\sum_{n\geq1}t^n\sum_{\lambda,\mu^0,...,\mu^k\vdash n}\tilde{C}^\lambda_{\mu^0,...,\mu^k}p_\lambda q^{(0)}_{\mu^0}...q^{(k)}_{\mu^k},$$
where $\tilde{C}^\lambda_{\mu^0,...,\mu^k}$ is the number of rooted connected orientable $k$-constellations,  with profile $(\lambda,\mu^0,...,\mu^k)$.
\end{thm}
\noindent In this paper, we give an analog result in the non-oriented case:

\begin{thm}\label{Thm b=1}
For every $k\geq1$, we have

$$
 (i)\hspace{1cm}  \tau^{(k)}_1(t,\mathbf{p},\mathbf{q}^{(0)},..,\mathbf{q}^{(k)})=1+\sum_{n\geq 1}t^n\sum_{\lambda,\mu^0,...,\mu^k\vdash n}\frac{|\mathfrak{F}^\lambda_{\mu^0,...,\mu^k}|}{2^{\ell(\lambda)}z_\lambda}p_\lambda q^{(0)}_{\mu^{0}}..q^{(k)}_{\mu^{k}},
$$
$$ (ii)\hspace{1cm} \Psi_1^{(k)}(t,\mathbf{p},\mathbf{q}^{(0)},..,\mathbf{q}^{(k)})=\sum_{n\geq1}t^n\sum_{\lambda,\mu^0,...,\mu^k\vdash n} C^\lambda_{\mu^0,...,\mu^k}p_\lambda q^{(0)}_{\mu^0}...q^{(k)}_{\mu^k},$$
where $\ell(\lambda)$ is the number of parts of $\lambda$ and $C^\lambda_{\mu^0,...,\mu^k}$ is the number of non-oriented rooted  connected $k$-constellations with profile $(\lambda,\mu^0,...,\mu^k)$.
\end{thm}

For $k=1$, these results have been proved in \cite{GJ96a}. In this case,  constellations are bipartite maps.
The generating series of non-oriented bipartite maps have been\footnote{Actually, the maps considered in \cite{GJ96a} are face colored maps called hypermaps. These maps are obtained from bipartite maps by duality.} investigated through the correspondence between bipartite maps and matchings on the one hand, and a relationship between matching enumeration and the structure coefficients of the double coset algebra of the Gelfand pair $(\mathfrak{S}_{2n},\mathfrak{B}_n)$ on the other hand. In this paper, we extend this development to constellations.

\subsection{Generalized Goulden and Jackson conjectures.}\label{sec Conj}
 We define the coefficients $c^\lambda_{\mu^0,...,\mu^k}$ and $h^\lambda_{\mu^0,...,\mu^k}$  for partitions $\lambda$,$\mu^0$,...,$\mu^k\vdash n\geq1$ such that 
\begin{equation}\label{defc}
    \tau^{(k)}_b(t,\mathbf{p},\mathbf{q}^{(0)},..,\mathbf{q}^{(k)})=1+\sum_{n\geq1}t^n\sum_{\lambda,\mu_0,...,\mu_k\vdash n}\frac{c^\lambda_{\mu_0,...,\mu_k}(b)}{z_\lambda(1+b)^{\ell(\lambda)}}p_\lambda q^{(0)}_{\mu_0}...q^{(k)}_{\mu_k},
\end{equation}
 and
\begin{equation}\label{defh}
    \Psi^{(k)}_b(t,\mathbf{p},\mathbf{q}^{(0)},..,\mathbf{q}^{(k)})=\sum_{n\geq1}t^n\sum_{\lambda,\mu_0,...,\mu_k\vdash n}h^\lambda_{\mu_0,...,\mu_k}(b)p_\lambda q^{(0)}_{\mu_0}...q^{(k)}_{\mu_k}.
\end{equation}
 For $k=1$, these coefficients were introduced by Goulden and Jackson. They have conjectured that these coefficients are non-negative integer polynomials in $b$ and that they enumerate respectively matchings and bipartite maps. Based on the particular cases of \cref{Thm b=0} and \cref{Thm b=1}, and on computational explorations we formulate a generalization of these conjectures for $k\geq1$. These generalized conjectures are somehow implicit in  \cite{CD20}.
First, we introduce a positivity conjecture.

\begin{conj}[Generalized positivity conjecture.]\label{Pos conj}
For every $k\geq1$ and partitions $\lambda,\mu^0,...,\mu^k$ of size $n\geq1$, the coefficients $c^\lambda_{\mu^0,...\mu^k}(b)$ and $h^\lambda_{\mu^0,...\mu^k}(b)$ are polynomials in $b$ with non-negative integer coefficients.
\end{conj}

It turns out that Goulden and Jackson's positivity conjecture for the coefficients $c$ (the case $k=1$) implies the conjecture for any $k\geq1$, since these coefficients satisfy a multiplicativity property (see \cref{prop mult}). Such property does not a priory exist for the coefficients $h^\lambda_{\mu^0,..,\mu^k}$, so the positivity for these coefficients is more general than the conjecture for $k=1$.  Using computer explorations, this conjecture has been tested  when $k\leq 5$ and $n\leq 12-k$. We attach\footnote{To access the attached files with acrobat reader, click on view/navigation panels/attachments. Alternatively, you can download the source files from arXiv.org}  the values of the coefficients $c$ and $h$ for $k\leq 5$ and $n\leq 9-k$. Given Theorems \ref{Thm b=0} and \ref{Thm b=1}, \Cref{Pos conj} is equivalent to the two following combinatorial conjectures: 
\begin{conj}[Generalized Matching-Jack conjecture]\label{match conj}
For every $k\geq1$ and partitions $\lambda,\mu^0,...,\mu^k$ of size $n\geq1$, there exists a function $\vartheta_\lambda$ on $\mathfrak{F}^\lambda_{\mu^0,...,\mu^k}$ with non-negative integer values, such that $\vartheta_\lambda(\delta_0,...,\delta_{k-1})$ is zero if and only if each one of the matchings $\delta_0$,...,$\delta_{k-1}$ is bipartite, and
$$c^\lambda_{\mu_0,...,\mu_k}(b)=\sum_{(\delta_0,...,\delta_{k-1})\in\mathfrak{F}^\lambda_{\mu^0,...,\mu^k}}b^{\vartheta_\lambda(\delta_0,...,\delta_{k-1})}.$$
\end{conj}

\begin{conj}[Generalized $b$-conjecture]\label{bconj}
For every $k\geq1$ and partitions $\lambda,\mu^0,...,\mu^k\vdash n\geq1$, there exists a function $\nu$ on connected rooted constellations with profile $(\lambda,\mu^0,...,\mu^k)$ with non-negative integer values, such that $\nu(\mathbf{M},c)$ is zero if and only if $(\mathbf{M},c)$ is orientable, and  $$h^\lambda_{\mu_0,...,\mu_k}(b)=\sum_{(\mathbf{M},c)}b^{\nu(\mathbf{M},c)},$$
where the sum runs over rooted connected $k$-constellations with profile $(\lambda,\mu^0,...,\mu^k)$.
\end{conj}

The conjectures introduced by Goulden and Jackson in \cite{GJ96b}
correspond to the case $k=1$ and are still open, despite many partial results \cite{La09,KV16,DF16,DF17,D17,KPV18,CD20}. Using basic facts about Jack symmetric functions, we can see that the quantities $c^\lambda_{\mu_0,...,\mu_k}(b)$ and $h^\lambda_{\mu_0,...,\mu_k}(b)$ are rational functions in $b$ with rational coefficients. The polynomiality of the quantities $c^\lambda_{\mu^0,...,\mu^k}$ and $h^\lambda_{\mu^0,...,\mu^k}$  for $k=1$ is not direct from the construction and has been open for twenty years. Do\l{}{\k{e}}ga and Féray have proved this polynomiality for coefficients $c^\lambda_{\mu,\nu}(b)$, see \cite{DF16}. Short after, they deduced the polynomiality of $h^\lambda_{\mu,\nu}(b)$, see \cite{DF17}. The proofs extend directly to any $k\geq1$.

\Cref{match conj} and \Cref{bconj} are closely related, since the functions $\tau^{(k)}$ and $\Psi^{(k)}$ are related by \cref{eqPsi}, and  $k$-constellations can be encoded by matchings, see \cref{prop1} and \cref{prop2}. However, we are not aware of any implication between them in the general case (even for $k=1$). This difficulty to pass from one conjecture to another is due to the fact that we divide by $z_\lambda(1+b)^{\ell(\lambda)}$ in the definition of $c^\lambda_{\mu^0,...,\mu^k}$ and we should take the logarithm and the derivative to pass from $\Psi^{(k)}$ to $\tau^{(k)}$. Combinatorially, this can be explained by the fact that constellations appearing in the sum of the $b$-conjecture are rooted while there is no natural way to "root" the elements of  the sets $\mathfrak{F}^\lambda_{\mu^0,...,\mu^k}$ appearing in the Matching-Jack conjecture. This notion of rooting will be discussed in more details in \cref{sec3}. Nevertheless, we were able to overcome these difficulties in the case of marginal sums and deduce a result for $c^\lambda_{\mu^0,...,\mu^k}$ from the analog result of \cite{CD20} on $h^\lambda_{\mu^0,...,\mu^k}$.

\subsection{The case of marginal sums.}
We consider the marginal sums of coefficients $c^\lambda_{\mu,\mu^1,..,\mu^k}$ and $h^\lambda_{\mu,\mu^1,..,\mu^k}$, defined as follows: for all $\lambda,\mu\vdash n\geq1$, and for all $l_1,...,l_k\geq1$ set 
\begin{equation}\label{eq marginal sums}
  h^\lambda_{\mu,l_1,..,l_k}:=\sum_{\mu^i\vdash n,\ell(\mu^i)=l_i}h^\lambda_{\mu,\mu^1,..,\mu^k} \hspace{1cm}\text{and}\hspace{1cm}c^\lambda_{\mu,l_1,..,l_k}:=\sum_{\mu^i\vdash n,\ell(\mu^i)=l_i}c^\lambda_{\mu,\mu^1,..,\mu^k},  
\end{equation}
where $\ell(\mu^i)$ denotes the number of parts of the partitions $\mu^i$.
The main motivation for formulating the generalized versions of Goulden and Jackson conjectures introduced above is the following theorem due to Chapuy and Do\l{}{\k{e}}ga \cite{CD20} that establishes the generalized version of $b$-conjecture for the marginal sums $h^\lambda_{\mu,l_1,..,l_k}$;
\begin{thm}[\cite{CD20}]\label{thm CD} 
For every $k\geq1$, partitions $\lambda,\mu\vdash n\geq1$ and $l_1,...l_k\geq1$, there exists a function $\nu$ on connected rooted constellations with profile $(\lambda,\mu^0,...,\mu^k)$ with non-negative integer values, such that
$\nu(\mathbf{M},c)$ is zero if and only if $(\mathbf{M},c)$ is orientable, and  $$h^\lambda_{\mu,l_1,...,l_k}=\sum_{(\mathbf{M},c)}b^{\nu(\mathbf{M},c)},$$
where the sum runs over rooted connected $k$-constellations with profile $(\lambda,\mu,\mu^1...,\mu^k)$ for some partitions $\mu^1,...,\mu^k$ satisfying $\ell(\mu^i)=l_i$.
\end{thm}
\noindent The second main result of this paper is an analog for the marginal sums $c^\lambda_{\mu,l_1,..,l_k}$:

\begin{thm}\label{Thm3}
For every $k\geq1$, partitions $\lambda,\mu\vdash n\geq1$, and $l_1,...l_k\geq1$, there exists a function $\vartheta_\lambda$ on $\mathfrak{F}^\lambda_{\mu^0,...,\mu^k}$ with non-negative integer values, such that $\vartheta_\lambda(\delta_0,...,\delta_{k-1})$ is zero if and only if each one of the matchings $\delta_0$,...,$\delta_{k-1}$ is bipartite, and  $$c^\lambda_{\mu,l_1,..,l_k}=\sum_{\ell(\mu^i)=l_i}\sum_{(\delta_0,...,\delta_{k-1})\in\mathfrak{F}^\lambda_{\mu,\mu^1...,\mu^k}}b^{\vartheta_\lambda(\delta_0,...,\delta_{k-1})}.$$
\end{thm}

As explained above, the implications between the $b$-conjecture and the Matching-Jack conjecture are still open problems. The proof that we give here to deduce \cref{Thm3} from \cref{thm CD} can not be applied in the general case of the conjectures, since it uses a property of symmetry of the statistic $\nu$ that appears in \cref{thm CD}, see \cref{eq U}. Note that for the other partial results established for these conjectures (the cases $b=0$, $b=1$, and polynomiality), we start by proving the result for the Matching-Jack conjecture and then we deduce it for the $b$-conjecture, this approach is reversed in the current case.

\subsection{Lassale's conjecture.}
For every partitions $\lambda\vdash n\geq1$, and $\mu\vdash m\geq1$ we define $\theta^{(\alpha)}_\mu(\lambda)$ as follows:
$$ \theta^{(\alpha)}_\mu(\lambda):=\left\{
    \begin{array}{ll}
       0, & \mbox{ if } n<m.\\
        \left[p_\mu\right]J^{(\alpha)}_\lambda, & \mbox{ if } n=m.\\
        \theta^{(\alpha)}_\mu(\lambda)=\binom{|\lambda|-|\mu|+m_1(\mu)}{m_1(\mu)}\theta^{(\alpha)}_{\mu,1^{n-m}}(\lambda) &\mbox{If $m<n$.}
    \end{array}
\right.$$
where $[p_\mu]$ is the extraction symbol with respect to the variable $\mathbf{p}$ and $m_1(\mu)$ is the number of parts equal to 1 in the partition $\mu$.

It is known that the coefficients $\theta_\mu^{(\alpha)}(\lambda)$ are polynomials in $\alpha$ with integer coefficients (see \cite{KS97} and the discussion of \cite[Section 3.6]{DF16}). Actually, it has been proved that these quantities are also polynomials in the multirectangular coordinates of the partition $\lambda$ (see \cite{Las08} for a precise definition). Lassale has conjectured that these polynomials satisfy some positivity properties;

\begin{conj}\label{conj Lassale}
For every partition $\mu$ of size $m$ such that $m_1(\mu)=0$, the quantities $(-1)^{|\mu|}z_\mu\theta_\mu^{(\alpha)}(\lambda)$ are polynomials in the parameters $(b,q_1,q_2,...,-r_1,-r_2,...)$ with non-negative integers coefficients, where $b:=\alpha-1$, and $(q_1,q_2..,r_1,r_2...)$ are the multirectangular coordinates of $\lambda$.
\end{conj}
In this paper, we consider the case where the partition $\lambda$ has a rectangular shape, \textit{i.e.} a partition with $q$ parts of size $r$, where $r,q\geq 1$. In this case, the multirectangular coordinates of $\lambda$ are given by $(q,r)$ and we write $\lambda=(q\times r)$, and $\theta^{(\alpha)}_\mu(q,r):=\theta^{(\alpha)}_\mu(\lambda)$. 

Using a recurrence formula for the coefficients $\theta^{(\alpha)}_\mu(q,r)$, Lassale has established in the same paper the positivity in \Cref{conj Lassale} for the rectangular case but not the integrality. In \cite{DFS14}, a combinatorial expression of these coefficients in terms of weighted bipartite maps was given.  However the weights considered are not integral.

\Cref{Thm Lassale conjecture} gives a complete answer to the rectangular case in \Cref{conj Lassale} by proving the integrality of the coefficients.  We obtain an expression of $(-1)^{|\mu|}z_\mu\theta^{(\alpha)}_\mu(q,-r)$ as a sum of bipartite maps with monomial weights in $q$, $-\alpha r$ and $b$. The approach we use here is different from the one used in \cite{Las08} and \cite{DFS14}. It is is based on the marginal case in the Matching-Jack conjecture (\Cref{Thm3}) and \Cref{lem recJack} that relates the Jack polynomials indexed by rectangular partitions to some specialisations of the function $\tau^{(1)}_b$.

\subsection{Outline of the paper.}
In \cref{Sec preliminaries} we introduce some necessary definitions and notation. In \cref{sec b=1}, we introduce a notion of labelling for non-oriented constellations, in order to construct a correspondence between $k$-constellations and matchings. Building on that, we prove \cref{Thm b=1}. \cref{sec3} is devoted to the proof of \cref{Thm3}. In \Cref{sec Lassale conjecture}, we prove the rectangular case in \Cref{conj Lassale}.
In \cref{sec Top degree}, we discuss some general properties of the generalized connection coefficients $c^\lambda_{\mu^0,...,\mu^k}$ and we give a new proof for the positivity of the top degree part of these coefficients.

\section{Preliminaries}\label{Sec preliminaries}

\subsection{Partitions.}\label{subsec Partitions}

A \textit{partition} $\lambda=[\lambda_1,...,\lambda_\ell]$ is a weakly decreasing sequence of integers $\lambda_1\geq...\geq\lambda_\ell>0$. The quantity $\ell$ is called the \textit{length} of $\lambda$ and is denoted $\ell(\lambda)$. The size of $\lambda$ is the integer $|\lambda|:=\lambda_1+\lambda_2+...+\lambda_\ell.$ If $n$ is the \textit{size} of $\lambda$, we say that $\lambda$ is a partition of $n$ and we write $\lambda\vdash n$. The integers $\lambda_1$,...,$\lambda_\ell$ are called the \textit{parts} of $\lambda$.  For every $i\geq1$, we denote by $m_i(\lambda)$ the number of parts of $\lambda$ which are equal to $i$. The partition $2\lambda$ is given by $2\lambda:=[2\lambda_1,2\lambda_2,...]$. 

We denote by $\mathcal{P}$ the set of all partitions, including the empty partition. For every partition $\lambda$ and $i\geq1$, we set $\lambda_i=0$ if $i > \ell(\lambda).$
The dominance partial ordering  $\leq$ on $\mathcal{P}$ is given by 
$$\mu\leq\lambda \iff |\mu|=|\lambda| \text{ and }\hspace{0.3cm} \mu_1+...+\mu_i\leq \lambda_1+...+\lambda_i \text{ for } i\geq1.$$

\noindent We identify a partition  $\lambda$ with its Young diagram defined by 
$$\lambda:=\{(i,j),1\leq i\leq \ell(\lambda),1\leq j\leq \lambda_i\}.$$
For every box $\Box:=(i,j)\in\lambda$, the  \textit{arm-length} of $\Box$ is given by $a_\lambda(\Box):=|\{(i,r)\in\lambda,r>j\}|$ and its \textit{leg-length} is given by $\ell_\lambda(\Box):=|\{(r,j)\in\lambda,r>i\}|$. Two $\alpha$-deformations of the hook-length product were introduced in \cite{stan89};
$$\hook_\lambda^{(\alpha)}:=\prod_{\Box\in\lambda}\left(\alpha a_\lambda(\Box)+\ell_\lambda(\Box)+1\right),\hspace{0.3cm}\hook_\lambda'^{(\alpha)}:=\prod_{\Box\in\lambda}\left(\alpha(a_\lambda(\Box)+1)+\ell_\lambda(\Box)\right).$$
With these notations, the classical hook-length product is given by (see e.g. \cite{stan89})
$$H_\lambda:=\hook_\lambda^{(1)}=\hook_\lambda'^{(1)}.$$
Finally, for every box $\Box:=(i,j)\in\lambda$, we define its $\alpha$-content by $c_\alpha(\Box):=\alpha(x-1)-(y-1)$.

\subsection{Matchings.} \label{subsec Matchings}
We introduce some notation related to matchings as defined in \cite{GJ96b}.
We recall that for every $n\geq1$, we set $\mathcal{A}_n :=\{1,\hat{1},...,n,\hat{n}\}$.
We also denote by $\mathfrak{F}_n$ the set of matchings on $\mathcal{A}_n$.
For all $\delta_1,...,\delta_r\in\mathfrak{F}_n$ we denote by $G(\delta_1,...,\delta_r)$ the multi-graph with vertex-set  $\mathcal{A}_n$, and edges all the pairs of $\delta_1\cup...\cup\delta_k$.
In the case $r=2$, we  note that all connected components of $G(\delta_1,\delta_2)$ are cycles of even size, so we can define $\Lambda(\delta_1,\delta_2)$ as the partition of $n$ obtained by taking half of the size of each connected component of $G(\delta_1,\delta_2)$.

\subsection{Symmetric functions and Jack polynomials.}\label{sec SymFun}
For the definitions and notation introduced in this subsection we refer to \cite{Mac95}.
We denote by $\mathcal{S}$ the algebra of symmetric functions with coefficients in $\mathbb Q$. For every partition $\lambda$, we denote $m_\lambda$ the monomial function, $p_\lambda$ the power-sum function and $s_\lambda$ the Schur function  associated to the partition $\lambda$.
If $\alpha$ is an indeterminate, let $\mathcal{S}_\alpha:=\mathbb{Q}[\alpha]\otimes\mathcal{S}$ the algebra of symmetric functions with rational coefficients in $\alpha$.
We denote by $\langle.,.\rangle_\alpha$ the $\alpha$-deformation of the Hall scalar product defined by 
$$\langle p_\lambda,p_\mu\rangle_\alpha=z_\lambda\alpha^{\ell(\lambda)}\delta_{\lambda,\mu},\text{ for }\lambda,\mu\in\mathcal{P}.$$ 
Macdonald \cite[Chapter VI.10]{Mac95} has proved that there exists a unique family of symmetric functions $(J_\lambda)$ in $\mathcal{S}_\alpha$ indexed by partitions, satisfying the following properties, called Jack polynomials;
$$\left\{ \begin{array}{ll}
       \text{Orthogonality: }
     &\langle J_\lambda,J_\mu\rangle_\alpha=0, \text{ for }\lambda\neq\mu,\\
      \text{Triangularity: }
     &[m_\mu]J_\lambda=0, \text{ unless }\mu\leq\lambda,\\
     \text{Normalization: }
     &[m_{1^n}]J_\lambda=n!, \text{ for }\lambda\vdash n,
\end{array}\right.$$
where $[m_\mu]J_\lambda$ denotes the coefficient of $m_\mu$ in $J_\lambda$, and $1^n$ is the partition with $n$ parts equal to~1.
For $\alpha=1$ and $\alpha=2$, the Jack polynomials are given by 
\begin{equation}
    J_\lambda^{(1)}=H_\lambda s_\lambda, \hspace{0.5cm } J_\lambda^{(2)}=Z_\lambda,
\end{equation}
where $Z_\lambda$ denotes the zonal polynomial associated to $\lambda$, see \cite[Chapters VI and VII]{Mac95}.
The squared norm of Jack polynomials can be expressed in terms of the deformed hook-length products, see \cite[Theorem 5.8]{stan89};

\begin{equation}\label{eq j alpha}
    j_\lambda^{(\alpha)}:=\langle J_\lambda,J_\lambda\rangle_\alpha=\hook^{(\alpha)}_\lambda\hook^{(\alpha)}_\lambda.
\end{equation}
In particular, we have
\begin{equation}\label{eq j}
    j_\lambda^{(1)}=H_\lambda^2 \hspace{0.3cm}\text{and}
    \hspace{0.3cm} j_\lambda^{(2)}=H_{2\lambda}.
\end{equation}

We conclude this subsection with the following theorem (see \cite[Equation 10.25]{Mac95});
\begin{thm}[\cite{Mac95}]\label{Jack formula}
For every $\lambda\in\mathcal{P}$, we have
$$J_\lambda(\underline{u})=\prod_{\Box\in\lambda}\left(u+c_\alpha(\Box)\right).$$
\end{thm}

\subsection{Maps.}\label{ssec Maps}
\begin{figure}[ht]
        \centering
    \includegraphics[width=.2\textwidth, center]{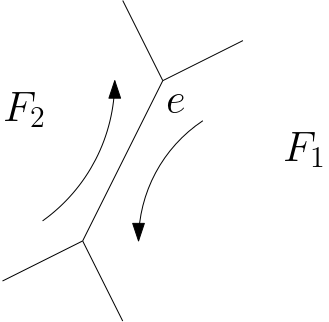}
    \caption{Consistent orientation from either side of an edge $e$.}
    \label{consistent Orientations}
   \end{figure}

We start by giving the definition of a map (see \cite[Definition 1.3.6]{LZ04}).
\begin{defi}
A \textit{connected map} is a connected graph embedded into a surface  such that all the connected components of the complement of the graph are simply connected. These connected components are called the \textit{faces} of the map. We consider maps up to homeomorphisms of the surface (see \cite[Definition 1.3.7]{LZ04}). A connected map is \textit{orientable} if it is the case for the underlying surface. A \textit{map} is an unordered collection of connected maps. A map is orientable if each one of its connected components is orientable. We will use the word \textit{non-oriented maps} for maps which are orientable or not.
\end{defi}
Another description of orientable maps is the following : a map is orientable if each one of its faces can be endowed with an orientation, such that for every edge $e$ of the map the two edge-sides of $e$ are oriented in opposite ways. In \cref{consistent Orientations} we have an edge $e$ whose sides are incident to two faces $F_1$ and $F_2$ (not necessarily distinct), and that are oriented in opposite ways. In this case we say that the orientation of the faces is \textit{consistent}.
In this paper, we will consider rooted maps, \textit{i.e.} maps with a distinguished oriented corner. 
We call  \textit{canonical orientation} of a rooted connected orientable map the unique orientation on the faces of the map which is consistent and such that the face containing the root is oriented by the root corner.

\subsection{\texorpdfstring{$k$}{}-Constellations}\label{sec cons}

In this subsection, we introduce the same notation related to constellations as in \cite{CD20}. 

\begin{figure}
\centering
\begin{minipage}{.53\textwidth}
  \centering
  \includegraphics[width=.69\linewidth]{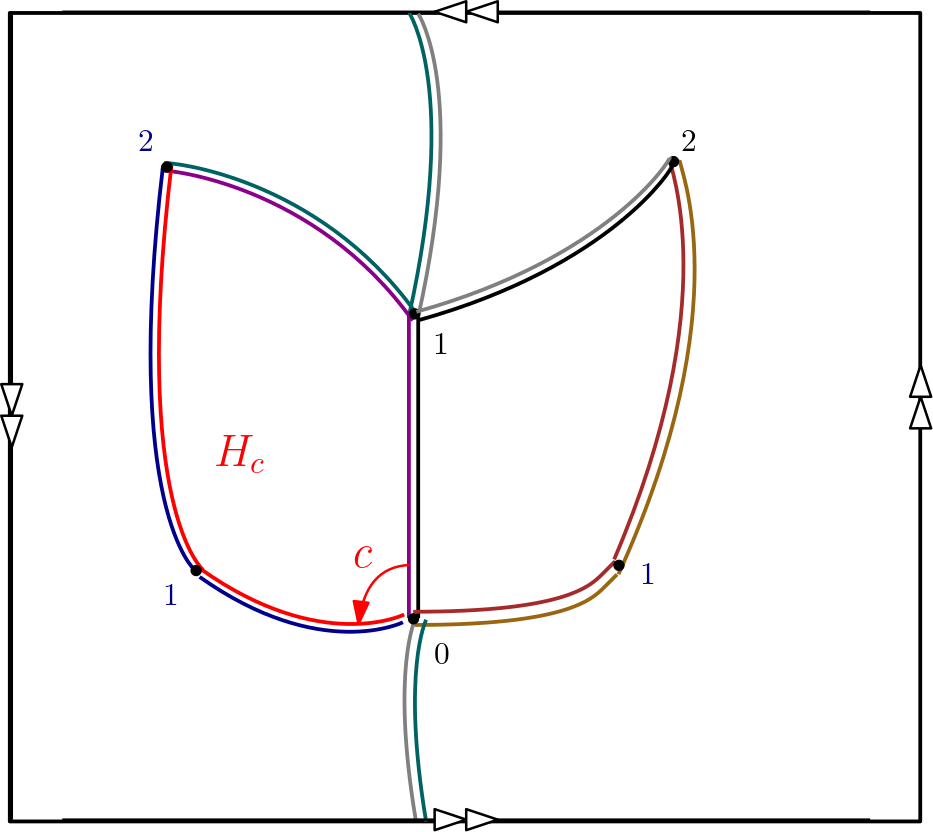}
   \captionof{figure}{Example of a rooted 2-constellation on the projective plane. $c$ is the root corner of the constellation, and $H_c$ is the root right-path. The 8 right-paths are represented in different colors. Note that the two sides of one edge are always in different right-paths.}
   \label{unlabelled 2 const}
\end{minipage}%
\begin{minipage}{.53\textwidth}
  \centering
   \includegraphics[width=.73\linewidth]{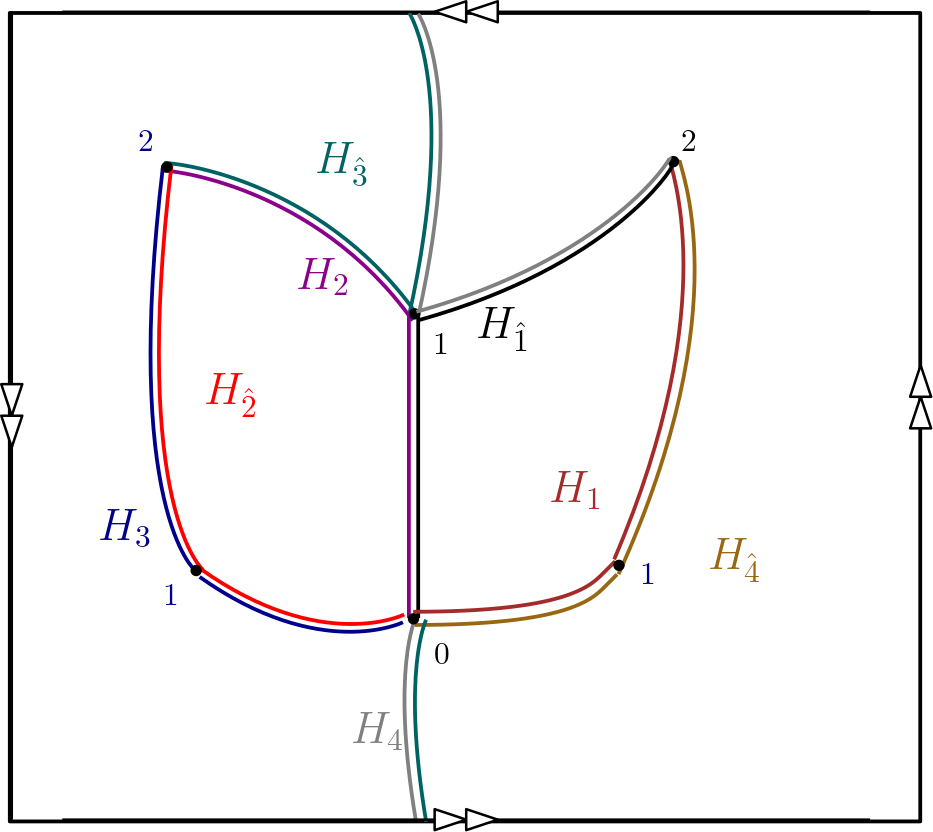}
  \captionof{figure}{An example of a labelling of the 2-constellation illustrated in \Cref{unlabelled 2 const}. The right-path labelled by $i$ is denoted $H_i$.}
  \label{Figure 1}
\end{minipage}
\end{figure}

We call \textit{a right-path} of a $k$-constellation $\mathbf{M}$, a path of length $k$ along the boundary of a face of $\mathbf{M}$ that separates a corner of color 0 with a consecutive corner of color $k$ incident to this face.
We recall that a connected $k$-constellation $\mathbf{M}$ is \textit{rooted} if it is equipped with a distinguished oriented corner $c$ of color 0. This oriented corner $c$ will be called the \textit{root} of the constellation.
This is equivalent to distinguish a right-path in $\mathbf{M}$, that will be the right-path following the root corner, see \Cref{unlabelled 2 const}. We will use the term root to design the root corner or the root right-path.
The rooted constellation will be denoted $(\mathbf{M},c)$.
We say that an edge is of color $\{i,i+1\}$ if its end points are of color $i$ and $i+1$. When $k=1$, $1$-constellations are  bipartite maps  and right-paths are edge-sides.

Since the number of right-paths contained in each face is even, we can define the degree of a face as half the number of its right-paths.
Similarly, we define  the degree of a vertex as half the number of right-paths that passes by this vertex (we can see that this is the number of edges incident to this vertex if it has color 0 or $k$, and half the number of edges incident to this vertex if it has color in $\{1,...,k-1\}$). 
We also define the size of a $k$-constellation $\mathbf{M}$ as half the number of its right-paths, it will be denoted $|\mathbf{M}|$. 
Therefore, for every $k$-constellation $\mathbf{M}$ and for every  color 
$0\leq i\leq k$, we have 
$|\mathbf{M}|=\sum_{v}\deg(v),$
where the sum runs over vertices of color $i$. We also have 
$|\mathbf{M}|=\sum_{f}\deg(f),$ where the sum runs over the faces of $\mathbf{M}$. We define the \textit{face-type} of a $k$-constellation as the partition obtained by reordering the degrees of the faces of $\mathbf{M}$. Similarly, for $i\in \{0,...,k\}$, the type of the vertices of color $i$, \textit{i.e.} the partition obtained by reordering the degrees of the vertices of color $i$.
We define the \textit{profile} of a $k$-constellation $\mathbf{M}$ as the $k+2$-tuple of partitions $(\lambda,\mu^0,...,\mu^k)$, such that $\lambda$ is the face-type of~$\mathbf{M}$, and for $i\in \{0,...,k\}$, $\mu^i$ the type of the vertices of color $i$.
If $\mathbf{M}$ is a $k$-constellation of size $n$, then $\lambda,\mu^0,...,\mu^k\vdash n$.

\begin{exe}
The 2-constellation of \Cref{unlabelled 2 const} has size $4$ and profile $([2,1,1],[4],[2,1,1],[2,2])$.
\end{exe}
Finally, we say that a $k$-constellation of size $n$ is labelled if it is equipped with a bijection between its right-paths and the set $\mathcal{A}_n=\{1,\hat{1},...,n,\hat{n}\}$. Labelled $k$-constellations will be denoted with a check : $\Check{\mathbf{M}}$.

\noindent We now compare the definition of constellations that we use here (given in \cref{def const}) to the usual description of orientable constellations given by hypermaps, (see e.g \cite{BMS00,Cha18,F16}). This correspondence between the two descriptions is mentioned in \cite{CD20} without details.

\textbf{Link with the usual definition of orientable constellations:}
We define a \textit{hypermap} as a map with faces colored in two colors, such that each edge separates two faces of different colors. The faces of one color are called the \textit{hyperedges} of the hypermap, and the faces of the other color are called the \textit{faces} of the hypermap.
The usual definition of orientable constellations is the following:
\begin{defi}[Usual definition of orientable constellations.]\label{def usual const} Let $k\geq 1$. An orientable $k$-constellation is an orientable hypermap with  vertices colored in the colors $\{0,1,2...,k\}$ with the following properties:
\begin{itemize}
    \item The degree of each \textit{hyperedge} is equal to $k+1$.
    \item The degree of each \textit{face} is a multiple of $k+1$.
    \item  There exits a consistent orientation of the faces, such that when we travel along a face in this orientation we read the colors $\{0,1,2...,k,0...\}$.
\end{itemize}
A connected orientable $k$-constellation is rooted if it has a distinguished hyperedge.
\end{defi}
Let us prove that this definition is equivalent to the definition of constellations that we use in this paper (see \Cref{def const}).

\begin{figure}[t]
\centering
\begin{subfigure}{.5\textwidth}
  \centering
  \includegraphics[width=.75\linewidth]{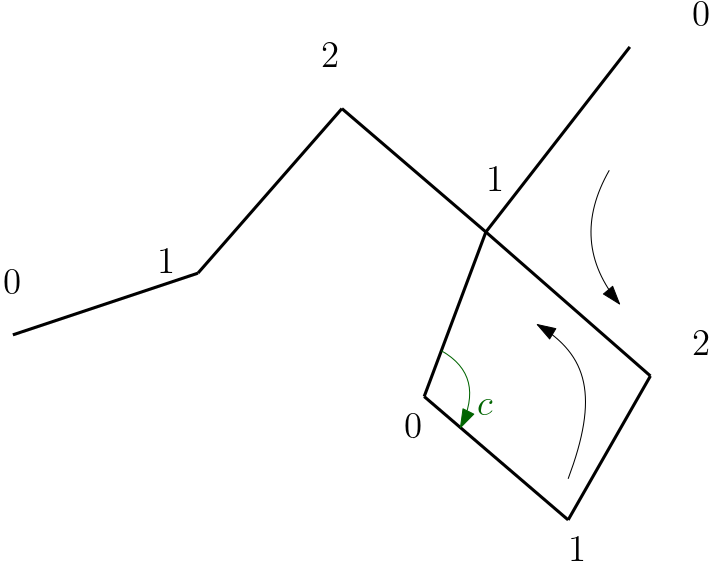}
  \label{fig:sub1}
\end{subfigure}%
\begin{subfigure}{.5\textwidth}
  \centering
  \includegraphics[width=.75\linewidth]{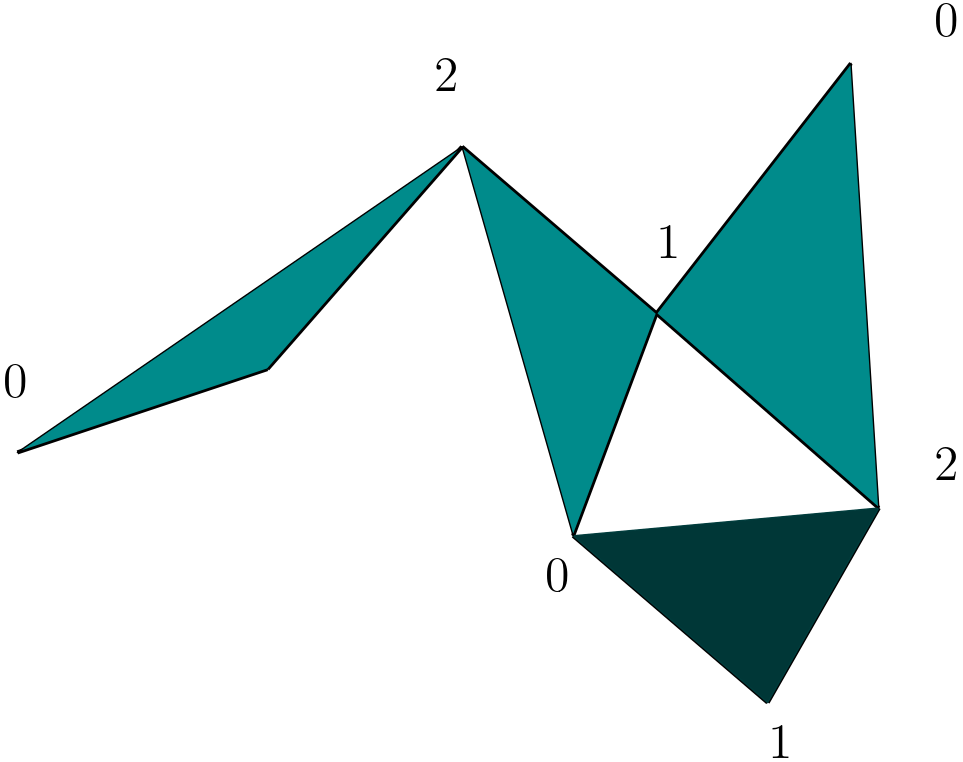}
  \label{fig:sub2}
\end{subfigure}
\caption{An example of an orientable 2-constellation in the plane, illustrated on the left with the description of \cref{def const}, and on the right with the description of hypermaps. The hyperedges are colored and the root hyperedge is represented with a darker color.}
\label{Orientable constellation}
\end{figure}

\noindent To each connected rooted orientable $k$-constellation (in the sense of \Cref{def const}), we can associate a hypermap as follows: we travel along each face with respect to the canonical orientation (see \Cref{ssec Maps}), and we add an edge between each corner of color 0 and the following corner of color $k$. In other terms, we close each right-path traversed from its corner of color 0 to its corner of color $k$ by adding an edge of color $(0,k)$, thus forming a face of degree $k$. Such face will be considered as a \textit{hyperedge} of the hypermap. The other faces of the map will be considered as \textit{faces}.
In \cref{Orientable constellation}, we have an example of this transformation illustrated on a planar 2-constellation. Let us prove that the map obtained is a hypermap; since the constellation is orientable, the orientations from either side of a given edge $e$ of the constellation are consistent. This implies that one of the two right-paths that contain $e$ is traversed form the corner of color 0 to the corner of color $k$ and the other right-path will be traversed in the opposite way. Only the first right-path will be transformed to a hyperedge. This proves that $e$ separates a \textit{hyperedge} and a \textit{face}. Hence the map obtained is a hypermap that satisfies the properties of \Cref{def usual const}.
Moreover, this constellation can be rooted by distinguishing the hyperedge associated to the root right-path. We thus recover the usual definition of orientable constellations.
Conversely, if we have a hypermap with the properties of \Cref{def usual const}, we can delete the edges of color $(0,k)$ to obtain a constellation as described in \cref{def const}.  
\begin{rmq}
Note that the orientability of the constellation is necessary to obtain a hypermap by the transformation described above. The description of orientable constellations with hypermaps has the advantage of being symmetric in the $k+1$ colors, while in the definition with right-paths the colors 0 and $k$ have a particular role. This lack of symmetry seems inevitable in the case of non-oriented constellations.
\end{rmq}

\section{The case \texorpdfstring{$b=1$}{}.}\label{sec b=1}

\subsection{Correspondence between constellations and tuples of matchings.}\label{sec correpondence}
The purpose of this subsection is to give a bijection between labelled $k$-constellations of size $n$ and $k+2$-tuples of matchings on $\mathcal{A}_n$. This is a generalization of the construction given in \cite{GJ96a} which corresponds to the case $k=1$.
\begin{defi}\label{Matchings}
If $\Check{\mathbf{M}}$ is a labelled $k$-constellation of size $n$, we define $\mathcal{M}(\mathbf{\Check{M}}):=(\delta_{-1},\delta_{0},...,\delta_k)$ as the $k+2$-tuple of matchings on  $\mathcal{A}_n$ defined as follows:
\begin{itemize}
    \item $\delta_{-1}$ (respectively $\delta_k$) is the matching whose pairs are the labels of right-paths of the same face, that have a corner of color $0$ (respectively $k$) in common.
\item For $i\in\{0,...,k-1\}$, $\delta_{i}$ is the matching whose pairs are the labels of right-paths having an edge of color $(i,i+1)$ in common.
\end{itemize}
\end{defi}

\noindent It is easy to see that the profile of a $k$-constellation $\mathbf{\Check{M}}$ can be determined by the associated matchings $\mathcal{M}(\mathbf{\Check{M}})$:
\begin{itemize}
    \item $\Lambda(\delta_{-1},\delta_k)$ is the face-type.
    \item For $i\in\{0,...,k\}$, $\Lambda(\delta_{i-1},\delta_i)$ is the type of vertices of color $i$.
\end{itemize}

Given a $k+2$-tuple of matchings $(\delta_{-1},..,\delta_k)$, we define its \textit{profile} as the $k+2$-tuple of partitions \\ $\Big(\Lambda(\delta_{-1},\delta_k),\Lambda(\delta_{-1},\delta_{0})...,\Lambda(\delta_{k-1},\delta_k)\Big)$.
We get from the previous remark  that $\Check{\mathbf{M}}$ and $\mathcal{M}(\Check{\mathbf{M}})$ have the same profile.

\begin{exe}\label{exe matchings}
 The labelled 2-constellation of \cref{Figure 1} is associated to the matchings $(\delta_{-1},\delta_0,\delta_1,\delta_2)$ below, $\color{red}\delta_{-1}$ in red, $\color{blue}\delta_0$ in blue, $\color{vert}\delta_1$ in green and $\delta_2$ in black.

\begin{center}
    
     \begin{tikzcd}
 4\arrow[rd,vert,dash]
&1 \arrow[ld, blue,dash]\arrow[ld,vert,dash,bend left =10]
& 2  \arrow[dl, blue,dash] \arrow[dr,vert,dash]
&3\arrow[d,dash]\arrow[r, bend left=30,red,dash] 
&4\arrow[dl, blue,dash]\\
\hat{4}
&\hat{1} \arrow[u,dash]\arrow[u,dash,red,bend right =10]
& \hat{2} \arrow[u,dash]\arrow[u,dash,red,bend right =10]\arrow[ur,vert,dash,bend right=10]\arrow[ur, blue,dash]
& \hat{3} \arrow[r, bend right=30,red,dash] 
& \hat{4} \arrow[u,dash]
\end{tikzcd}\\
\end{center}
\end{exe}

Conversely, if $(\delta_{-1},...,\delta_k)$ is a $k+2$-tuple of matchings of $\mathcal{A}_n$ we can construct a labelled $k$-constellation $\Check{\mathbf{M}}$ such that $\mathcal M({\Check{\mathbf{M}}}):=(\delta_{-1},...,\delta_k)$ (this construction is described with more details in \cite[Section 3.1]{DFS14} in the case $k=1$): 
\begin{itemize}
    \item For each connected component $C$ of the graph $G(\delta_{-1},\delta_k)$ of size $2r$ we consider a polygon consisting of $2r$ right-paths labelled by vertices in $C$, as follows: two right-paths have a vertex of color 0 in common (respectively of color $k$) if and only if their labels in $G(\delta_{-1},\delta_k)$ are connected by $\delta_{-1}$ (respectively $\delta_k$).
    \item For each $0\leq i\leq k-1$, and for each edge $e=(j,k)$ of the matching $\delta_i$ (where $j,k\in\mathcal{A}_n$), we glue the two edge-sides of color $(i,i+1)$ of the right-paths labelled by $j$ and $k$.
    \end{itemize}

    \begin{exe}
    The construction of the 2-constellation of  \Cref{Figure 1} by gluing polygons with respect to the matchings of \cref{exe matchings} is illustrated in \cref{Polygons gluing}.
    
    \begin{figure}[ht]
   \includegraphics[width=.5\textwidth, center]{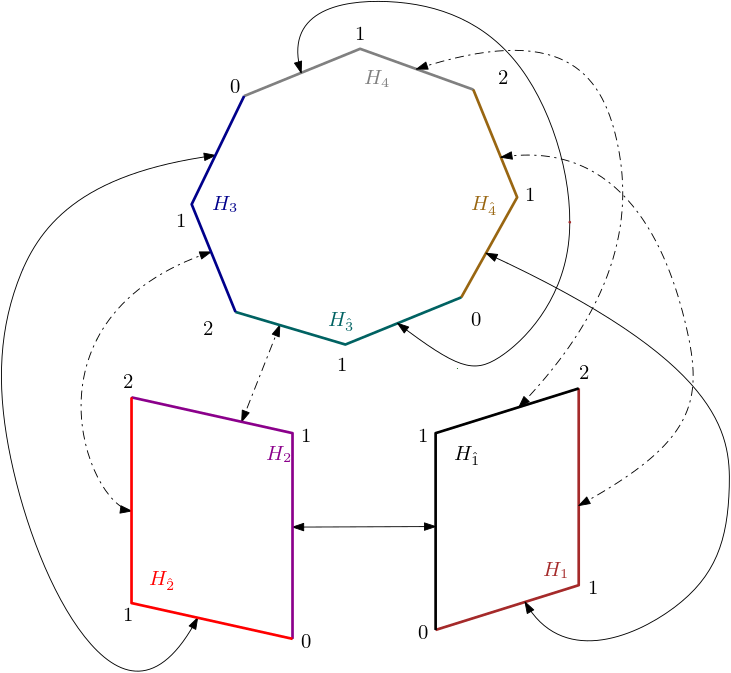}
   \caption{Polygons obtained from matchings $\delta_{-1}$ and $\delta_{2}$ of \cref{exe matchings}. Continuous arrows illustrate how to glue edge-sides of color $(0,1)$ with respect to matchings $\delta_0$ and dotted arrows illustrate how to glue edge-sides of color $(1,2)$ with respect to $\delta_1$.}
   \label{Polygons gluing}
   \end{figure}

    \end{exe}
    
\noindent 
From the definition of a map, we know that the faces of a constellation are isomorphic to open polygons. This implies every map can be obtained by gluing polygons as above (see \cite[Construction 1.3.20]{LZ04} for a complete proof in the orientable case). We deduce the following proposition:
 \begin{thm}\label{prop1}
 For $\lambda,\mu^0,...,\mu^k\vdash n$, the map $\Check{\mathbf{M}}\longmapsto \mathcal{M(\mathbf{\Check{M}})}$ is a bijection between labelled $k$-constellation with profile $(\lambda,\mu^0,...,\mu^k)$ and $k+2$-tuples of matchings on $\mathcal{A}_n$ with the same profile.
 \end{thm}

\noindent Finally, we use the previous correspondence between constellations and matchings to introduce the notion of duality that will be useful in \cref{sec3}.
\begin{defi}\label{def duality}
Let $(\mathbf{M},c)$ be a rooted $k$-constellation. We define the dual constellation $(\tilde{\mathbf{M}},\tilde c)$ as follows. First, we choose a labelling of $\mathbf{M}$ such that the root right-path is labelled by 1, we obtain a labelled constellation $\Check{\mathbf{M}}$. Let $(\delta_{-1},\delta_0,...,\delta_{k}):=\mathcal{M}(\Check{\mathbf{M}})$. Then, We define the labelled constellation $\Check{\mathbf{M}}'$ such that $(\delta_{-1},\delta_k,...,\delta_{0})=\mathcal{M}(\Check{\mathbf{M}}')$ (\textit{i.e.} we exchange the matchings $\delta_i\leftrightarrow \delta_{k-i}$ for $0\leq i\leq k$).  Finally we forget the labels of $\Check{\mathbf{M}}'$ except for the label 1. We obtain a rooted constellation $(\tilde{ \mathbf{M}},\tilde c)$. It is clear that $(\tilde{ \mathbf{M}},\tilde c)$ does not depend on the labelling chosen for $(\mathbf{M},c)$.
\end{defi}

\noindent One can check that this definition is consistent with the definition of duality given in \cite[Definition 2.4]{CD20}. 
\begin{rmq}\label{rmq duality}
It is straightforward from the definition that duality is an involution that exchanges faces with vertices of color 0, and vertices of color $i$ with vertices of color $k+1-i$, for $1\leq i \leq k$. More precisely, for every partitions $\lambda,\mu^0,...,\mu^k$ duality is a bijection between $k$-constellations with profile $(\lambda,\mu^0...,\mu^k)$ and constellations with profile $(\mu^0,\lambda,\mu^k,...,\mu^1)$.
\end{rmq}
\begin{rmq}
It is also possible, using matchings, to generalize this notion of duality in order to exchange colors in all possible ways, while controlling the profile as in the previous remark. However, these generalizations do not have a simple description in terms of maps.
\end{rmq}

\subsection{The Gelfand pair \texorpdfstring{$(\mathfrak S_{2n},\mathfrak{B}_n$)}{}}
In this subsection, we give some results that will be useful in the proof of \cref{Thm b=1} item \textit{(i)}. We follow the computations given in \cite{GJ96a} when $k=1$, we recall the most important steps of this proof and give a generalized version for the key lemmas. For this purpose we need to recall some results on the Gelfand pair  $(\mathfrak S_{2n},\mathfrak{B}_n)$ (see \cite[Section VII.2]{Mac95}).
We consider $\mathfrak{S}_{2n}$ as the permutation group of the set $\mathcal{A}_n:=\{1,\hat{1},...,n,\hat{n}\}$. We define the following action of $\mathfrak{S}_{2n}$ on $\mathfrak{F}_n$, the set of matchings on $\mathcal{A}_n$.
\begin{defi}
Let $\sigma\in\mathfrak{S}_{2n}$ and $\delta\in\mathfrak{F}_n$. We define  $\sigma.\delta$ as the matching of $\mathfrak{F}_n$ such that 
$(i,j)$ is a pair of $\sigma.\delta$ if and only if $(\sigma^{-1}i,\sigma^{-1}j)$ is a pair of $\delta$.
\end{defi}   
\noindent This action is both transitive and faithful.
We define the hyperoctahedral group $\mathfrak{B}_n$ as the stabilizer subgroup of the matching $\varepsilon$. One has that $|\mathfrak{B}_n|=n!2^n$.

\begin{defi}
Let $\sigma\in\mathfrak{S}_{2n}$. We define the coset-type of $\sigma$ as the partition of $n$ given by $\Lambda(\varepsilon,\sigma\varepsilon)$.
\end{defi}

The double cosets $\mathfrak{B}_n\backslash\mathfrak{S}_{2n}/\mathfrak{B}_n$ can be indexed by the partitions of $n$, in fact for all $\sigma,\tau\in\mathfrak{S}_{2n}$, one has
$\mathfrak{B}_n\sigma\mathfrak{B}_n=\mathfrak{B}_n\tau\mathfrak{B}_n$ if and only if $\sigma$ and $\tau$ have the same coset-type (see \cite[Lemma 3.1]{HSS92}). We denote $\mathcal{K}_\lambda$ the class of $\mathfrak{B}_n\backslash\mathfrak{S}_{2n}/\mathfrak{B}_n$ indexed by the partition $\lambda$, \textit{i.e.} the class of permutations of coset-type $\lambda$ and $K_\lambda\in\mathbb{C}\mathfrak{S}_{2n}$ defined by $$K_\lambda:=\sum_{\sigma\in\mathcal{K}_\lambda}\sigma.$$
These sums are the basis of a commutative subalgebra of $\mathbb{C}\mathfrak S_{2n}$, the Hecke algebra of the Gelfand pair $(\mathfrak S_{2n},\mathfrak{B}_n)$ (see \cite[Section VII.2]{Mac95}).
 Hence for $\lambda,\mu^0,...,\mu^k\vdash n$, we define $a^{\lambda}_{\mu^0,\mu^1,..\mu^k}$ such that 
 $$\prod\limits_{0\leq i\leq k} K_{\mu^i}=\sum_{\lambda\vdash n} a^{\lambda}_{\mu^0,\mu^1,..\mu^k}K_{\lambda}.$$
For $\sigma \in\mathcal{K}_\lambda$, $a^\lambda_{\mu^0,...,\mu^k}$ can be interpreted as the number of factorizations $\sigma=\sigma_0...\sigma_k$ where $(\sigma_0,...,\sigma_k)\in\mathcal{K}_{\mu^0}\times...\times\mathcal{K}_{\mu^k}$. We deduce that 
 \begin{equation}
     a^{\lambda}_{\mu^0,\mu^1,..\mu^k}|\mathcal K_{\lambda}|=\big|\left\{(\sigma_0,...,\sigma_k)\in\mathcal{K}_{\mu^0}\times...\times\mathcal{K}_{\mu^k} \text{ such that }\sigma_0...\sigma_k\in\mathcal{K}_\lambda\right\}\big|.
 \end{equation}
 For every $\lambda\vdash n$, there exist $\frac{n!}{z_\lambda}2^{n-\ell(\lambda)}$ matchings $\delta$ such that $\Lambda(\varepsilon,\delta)=\lambda$, see \cite[Proposition 5.2]{GJ96b}. On the other hand, using the fact that the action of $\mathfrak{S}_{2n}$ on $\mathfrak{F}_n$ is transitive, we can see that
 \begin{equation}\label{eq action}
     |\left\{\sigma\in\mathfrak{S}_{2n}
    \text{ such that } \sigma.\varepsilon=\delta_\lambda\right\}|=|\mathfrak{B}_n|.
 \end{equation}
 We deduce that 
 \begin{equation}\label{eq Kcal}
     |\mathcal{K}_\lambda|=|\mathfrak{B}_n|\frac{n!}{z_\lambda}2^{n-\ell(\lambda)}=\frac{|\mathfrak{B}_n|^2}{z_\lambda2^{\ell(\lambda)}}.
 \end{equation}
The coefficients $a^{\lambda}_{\mu^0,\mu^1,..\mu^k}$ are related to the size of the sets $\mathfrak{F}^\lambda_{\mu^0,...,\mu^k}$, defined in \cref{mathfrakF}. The following lemma is a generalization of \cite[Lemma3.2]{HSS92}.

 \begin{lem}\label{lem2}
 For $\lambda,\mu^0,...,\mu^k\vdash n$, we have
 $$|\mathfrak F^{\lambda}_{\mu^0,...,\mu^k}|=\frac{a^{\lambda}_{\mu^0,\mu^1,..\mu^k}}{|\mathfrak{B}_n|^k}.$$
 \begin{proof}
 
We define $\mathcal{E}\subset \mathfrak{S}_{2n}^{k+1}$, by $\mathcal{E}:=\left\{(\sigma_0,...,\sigma_k)\in\mathcal{K}_{\mu^0}\times...\times\mathcal{K}_{\mu^k} \text{ such that } \sigma_0...\sigma_k.\varepsilon=\delta_\lambda\right\}$,\\
where $\delta_\lambda$ is the matching defined in \cref{ssMatchings}. 
For every permutation $\sigma_\lambda$ such that $\sigma_\lambda.\varepsilon=\delta_\lambda$, one has that $\sigma_\lambda\in\mathcal{K}_\lambda$, and that $\sigma_\lambda$ has $a^\lambda_{\mu^0,...,\mu^k}$ factorizations of the form $\sigma_\lambda=\sigma_0...\sigma_k$ where $(\sigma_0,...,\sigma_k)\in\mathcal{K}_{\mu^0}\times...\times\mathcal{K}_{\mu^k}$. Using \cref{eq action}, we get 
 \begin{equation}\label{1}
 |\mathcal{E}|=|\mathfrak{B}_n|a^{\lambda}_{\mu^0,\mu^1,..\mu^k}.
 \end{equation}
 We now consider the map 
 \begin{align*}
 \psi:\hspace{1cm}
    &\mathcal E\hspace{1,3cm}\longrightarrow \mathfrak{F}_n^{k} \\
    &(\sigma_0,...,\sigma_k)\longmapsto (\sigma_0.\varepsilon,\sigma_0\sigma_1.\varepsilon,...,\sigma_0...\sigma_{k-1}.\varepsilon).
 \end{align*}
 For all $i\in\{0,...,k\}$, we have $$\Lambda(\sigma_0\sigma_{1}...\sigma_{i}.\varepsilon,\sigma_{0}...\sigma_{i-1}.\varepsilon)=\Lambda(\sigma_{i}.\varepsilon,\varepsilon)=\mu^i,$$ since $\sigma_i\in\mathcal{K}_{\mu^i}$.
 Hence, $\psi(\mathcal{E})\subseteq\mathfrak F^{\lambda}_{\mu^0,...,\mu^k}$. Let $(\delta_0,..,\delta_{k-1})\in\mathfrak F^{\lambda}_{\mu^0,...,\mu^k}$. There exists $(\sigma_0,...,\sigma_k)$ such that $\sigma_0.\varepsilon=\delta_0$, $\sigma_k.\delta_{k-1}=\delta_\lambda$ and for $i \in\{1,...,k-1\}$, $\sigma_i. \delta_{i-1}=\delta_i$. Then   
 $(\sigma_0,...,\sigma_k)\in\mathcal{E}$ and $\psi(\sigma_0,...,\sigma_k)=(\delta_0,...,\delta_{k-1})$, proving that  $\psi(\mathcal{E})=\mathfrak F^{\lambda}_{\mu^0,...,\mu^k}$. Moreover, $\psi(\sigma_0,...,\sigma_k)=\psi(\sigma'_0,...,\sigma'_k)$ if and only if there exist $\tau_0,...,\tau_k\in\mathfrak{B}_n$ such that for $i\in\{0,...,k\}$ 
 $$\sigma'_0...\sigma'_i=\sigma_0\sigma_1...\sigma_i\tau_i.$$
We deduce that  
\begin{equation}\label{2}
 |\psi^{-1}(\psi(\sigma_0,...,\sigma_k))|=|\mathfrak{B}_n|^{k+1}.
 \end{equation}
 Equations \eqref{1} and \eqref{2} conclude the proof.
 \end{proof}
 \end{lem}

We shall now establish the connection between the coefficients $a^\lambda_{\mu^0,...,\mu^k}$ and the Jack polynomials for $\alpha=2$. 
To this purpose we define for all $\lambda,\nu\vdash n$: 
$$\phi^\nu(\lambda):=\sum_{\sigma \in \mathcal K_\lambda}\chi^{2\nu}(\sigma),$$
where $\chi^{2\nu}$ is the irreducible character of the symmetric group  indexed by the partition $2\nu$. We also introduce the orthogonal idempotents of the Hecke algebra of $(\mathfrak S_{2n},\mathfrak{B}_n)$ that can be defined as follows (see \cite[Eq. (3.5)]{HSS92});

\begin{equation} \label{E}
E_\nu=\frac{1}{H_{2\nu}}\sum_{\lambda\vdash n} \frac{1}{|\mathcal K_\lambda|}\phi^\nu(\lambda)K_\lambda,
\end{equation}
They satisfy the property $E_\nu E_\rho=\delta_{\nu\rho}E_\nu$ for each $\nu,\rho\vdash n$, where $\delta$ is the Kronecker delta. \cref{E} can be inverted as follows (see \cite[Eq. (3.3)]{HSS92}):
\begin{equation} \label{K}
K_\lambda=\sum_{\nu\vdash n} \phi^\nu(\lambda)E_\nu.
\end{equation}
The following lemma is a generalization of \cite[Lemma 3.3]{HSS92}.
\begin{lem}\label{lem coefb}
For each partitions $\lambda,\mu^0,...,\mu^k\vdash n\geq1$, we have

 \begin{equation*}
     a^{\lambda}_{\mu^0,\mu^1,..\mu^k}=\frac{1}{|\mathcal{K}_{\lambda}|}\sum_{\nu\vdash n}\frac{1}{H_{2\nu}}\phi^\nu(\lambda)\phi^\nu(\mu^0)...\phi^\nu(\mu^k).
 \end{equation*}
 \begin{proof}
 Using \cref{K} we can write
 $$\prod\limits_{0\leq i\leq k}K_{\mu^i}=\sum_{\nu\vdash n}\phi^{\nu}(\mu^0)...\phi^\nu(\mu^k)E_\nu.$$
 We use \cref{E} to extract the coefficient of $K_\lambda$ from the last equality to obtain $a^{\lambda}_{\mu^0,\mu^1,..\mu^k}$.
 \end{proof}
\end{lem}

When $\alpha=2$, the Jack polynomials are called zonal polynomials and denoted by $Z_\theta$, see \cite[Chapter VII]{Mac95}. They can be expressed in the basis of power-sum functions as follows ; for every $\theta\vdash n$ one has
\begin{equation}\label{Zon}
  Z_\theta=\frac{1}{|\mathfrak B_n|}\sum_{\mu\vdash n} \phi^\theta(\mu)p_\mu.
\end{equation}
We are now ready to prove \cref{Thm b=1}.

\begin{proof}[Proof of \cref{Thm b=1} item.(i)]
 
   For $\alpha=2$, the function $\tau_1^{(k)}$ has the following expression; see Equations \eqref{eqtau} and \eqref{eq j}.
   $$\tau^{(k)}_1(t,\mathbf{p},\mathbf{q}^{(0)},..,\mathbf{q}^{(k)})=\sum_{n\geq0}t^n\sum_{\theta\vdash n}\frac{1}{H_{2\theta}}Z_\theta(\mathbf{p}) Z_\theta(\mathbf{q}^{(0)})...Z_\theta(\mathbf{q}^{(k)}).$$
   
   Using \cref{Zon} and \cref{lem coefb}, this can be rewritten as
   \begin{align*}
   \tau^{(k)}_1(t,\mathbf{p},\mathbf{q}^{(0)},..,\mathbf{q}^{(k)})
   &=\sum_{n\geq0}t^n\sum_{\theta\vdash n}\frac{1}{H_{2\theta}|\mathfrak{B}_n|^{k+2}}\sum_{\lambda,\mu^0,...,\mu^k\vdash n}\phi^\theta(\lambda)p_\lambda\phi^\theta(\mu^0)q^{(0)}_{\mu^0}...\phi^\theta(\mu^k)q^{(k)}_{\mu^k}\\
   &=\sum_{n\geq0}t^n\sum_{\lambda,\mu^{0},..\mu^{k}\vdash n}a^{\lambda}_{\mu^0,\mu^1,..\mu^k}\frac{|\mathcal{K_\lambda|}}{|\mathfrak{B}_n|^{k+2}}p_{\lambda}q^{(0)}_{\mu^{0}}..q^{(k)}_{\mu^{k}}.
    \end{align*} 
Finally, we use \cref{lem2} and \cref{eq Kcal} to conclude.
\end{proof}
\noindent Before deducing item (\textit{ii}) of \cref{Thm b=1}, we introduce the following notation; if $\mathbf{M}$ is a a $k$-constellation  with profile $(\lambda,\mu^0,...,\mu^k)$, we define the \textit{marking}\footnote{What is called marking in \cite{CD20} will be called marginal marking in this paper, see \Cref{ssec generalities}.} of $\mathbf{M}$ as the monomial
$$\tilde \kappa(\mathbf{M}):=p_{\lambda} q^{(0)}_{\mu^0} q^{(1)}_{\mu^1}...q^{(k)}_{\mu^k}.$$
We define the marking of a labelled constellation as the marking of the underlying constellation.
\noindent \cref{Thm b=1} item (\textit{ii}) can be reformulated as follows :

$$\Psi^{(k)}_1(t,\mathbf{p},\mathbf{q}^{(0)},...,\mathbf{q}^{(k)})=\sum_{(\mathbf{M},c)}t^{|\mathbf{M}|} \tilde {\mathbf{\kappa}}(\mathbf{M}),$$
where the sum runs over non-oriented rooted connected constellations.
\begin{proof}[Proof of \cref{Thm b=1} item(ii)]
\cref{Thm b=1} item.(\textit{i}) can be rewritten as follows; 
\begin{align*}
    \tau^{(k)}_1(t,\mathbf{p},\mathbf{q}^{(0)},..,\mathbf{q}^{(k)})&=\sum_{n\geq0}\frac{t^n}{(2n)!}\sum_{\lambda,\mu^{0},..\mu^{k}\vdash n}|\mathfrak F^{\lambda}_{\mu^0,..,\mu^k}|\frac{(2n)!}{n!2^n}2^{n-l(\lambda)}\frac{n!}{z_\lambda}p_{\lambda}q^{(0)}_{\mu^{0}}..q^{(k)}_{\mu^{k}},
\end{align*}
 On the other hand, the number of $k+2$-tuple of matchings $(\delta_{-1},...,\delta_k)$ with profile  $(\lambda,\mu^0...,\mu^{k})$ is given by  $\frac{(2n)!}{n!2^n}2^{n-l(\lambda)}\frac{n!}{z_\lambda}|\mathfrak F^{\lambda}_{\mu^0,..,\mu^k}|$; we have $\frac{(2n)!}{n!2^n}$ choices for $\delta_{-1}$, $\frac{n!}{z_\lambda}2^{n-l(\lambda)}$ choices for $\delta_k$ and $|\mathfrak F^{\lambda}_{\mu^0,..,\mu^k}|$ choices for the other matchings.
Using the description of labelled $k$-constellations with matchings (see \cref{prop1}) we obtain 
$$\tau^{(k)}_1(t,\mathbf{p},\mathbf{q}^{(0)},..,\mathbf{q}^{(k)})=\sum_{\Check{\mathbf{M}}}\frac{t^{|\Check{\mathbf{M}}|}}{(2|\Check{\mathbf{M}}|)!} \tilde {\mathbf{\kappa}}(\Check{\mathbf{M}}),$$
where the sum is taken over labelled $k$-constellations, connected or not.
Since the marking $\tilde{\kappa}(\mathbf{M})$ is multiplicative on the connected components of $\mathbf{M}$, we can apply the logarithm on the last equality in order to obtain the exponential generating series of connected labelled constellations (we use here the exponential formula for labelled combinatorial classes see e.g. \cite[Chapter II]{FS09}). When we forget all the labels in a connected rooted constellation except for the label "1", we obtain a constellation with a marked right-path that we can consider as a rooted constellation, see \cref{def const}. As each rooted constellation of size $n$ can be labelled in $(2n-1)!$ ways, we obtain 
$$\log\left(\tau^{(k)}_1(t,\mathbf{p},\mathbf{q}^{(0)},..,\mathbf{q}^{(k)})\right)=\sum_{(\mathbf{M},c)}\frac{t^{|\mathbf{M}|}}{2|\mathbf{M}|} \tilde {\mathbf{\kappa}}(\mathbf{M}),$$
where the sum runs over connected rooted constellations.
We conclude the proof by applying  $2\frac{t\partial}{\partial t}$ on the last equality.
\end{proof}

\section{Matching-Jack conjecture for marginal sums}\label{sec3}

\noindent The purpose of this section is to give a proof for \cref{Thm3}.
\subsection{Notation}\label{ssec generalities}
We consider two sequences of variables $\mathbf{p}=(p_1,p_2,...)$, $\mathbf{q}=(q_1,q_2,...)$ and $k$ variables $u_1$,...$u_k$. For a variable $u$ we denote   $\underline{u}:=(u,u,...)$. From the definition of the marginal sums $c^\lambda_{\mu,l_1...,l_k}$ and $h^\lambda_{\mu,l_1...,l_k}$ (see \Cref{eq marginal sums}), we have
\begin{equation*}\label{defc'}
    \tau_b^{(k)}(t,\textbf{p},\textbf{q},\underline{u_1},\underline{u_2},...,\underline{u_k})=\sum_{n\geq1}t^n\sum_{\lambda,\mu\vdash n}\sum_{l_1,...,l_k\geq 1}\frac{c^\lambda_{\mu,l_1...,l_k}(b)}{z_\lambda(1+b)^{\ell(\lambda)}}p_\lambda q_{\mu}u_1^{l_1}...u_k^{l_k},
\end{equation*}
\begin{equation*}
    \Psi_b^{(k)}(t,\textbf{p},\textbf{q},\underline{u_1},\underline{u_2},...,\underline{u_k})=\sum_{n\geq1}t^n\sum_{\lambda,\mu\vdash n}\sum_{l_1,...,l_k\geq 1}h^\lambda_{\mu,l_1...,l_k}(b)p_\lambda q_{\mu}u_1^{l_1}...u_k^{l_k}.
\end{equation*}
If $\mathbf{M}$ is a $k$-constellation with profile $(\lambda,\mu^0,\mu^1,...,\mu^k)$, we define the \textit{marginal marking} of $\mathbf{M}$ by $$\kappa(\mathbf{M}):=p_{\lambda} q_{\mu^0}u_1^{\ell(\mu^1)}...u_k^{\ell(\mu^k)},$$
and we say that $\big(\lambda,\mu^0,\ell(\mu^1),...,\ell(\mu^k)\big)$ is \textit{the marginal profile} of $\mathbf{M}$.  We can formulate \Cref{thm CD} as follows
\begin{thm}[\cite{CD20}]\label{CD} 
For every $k\geq1$, we have 
$$\Psi_b^{(k)}(t,\textbf{p},\textbf{q},\underline{u_1},\underline{u_2},...,\underline{u_k})=\sum_{(\textbf{M},c)} \kappa(\textbf{M}) t^{|\textbf{M}|}b^{\nu(\textbf{M},c)},$$
where the sum is taken over rooted connected $k$-constellations and $\nu(\mathbf{M,c})$ is a non-negative integer which is zero if and only if $(\mathbf{M},c)$ is orientable.
\end{thm}

\begin{defi}
For a class of vertex-colored maps, we call a $b$-weight a function $\rho$ that has values in 
$\mathbb{Q}[b]$ which has the two following properties:
\begin{itemize}
    \item $\rho(\mathbf{M})=1$ if and only if $\mathbf{M}$ is orientable.
    \item When we take $b=1$ we have $\rho(\mathbf{M})=1$.
\end{itemize}
  Moreover, we say that a $b$-weight  $\rho$ is integral if for every map $\mathbf{M}$ one has that $\rho(\mathbf{M})$ is a monomial in $b$. 
\end{defi}
\noindent With the definition above, the quantity $b^{\nu(\mathbf{M},c)}$ that appears in \cref{CD} is an integral $b$-weight on connected rooted constellations.
In \cref{sec bweights face-labelled}, we will consider $b$-weights on face-labelled constellations. 
\begin{rmq}\label{rmq thmCD20}
There is not a unique $b$-weight satisfying \Cref{thm CD}, see \cite[Theorem 5.10]{CD20}. In particular there exist non integral $b$-weights with this property. In this section, we fix once and for all an integral $b$-weight $b^{\nu(\mathbf{M},c)}$. 
\end{rmq}
For every $\lambda,\mu\vdash n$ and $l_1,...,l_k\geq1$ we define 
$$\mathfrak F^\lambda_{\mu,l_1,...,l_k}:=\bigcup\limits_{\mu^i\vdash n,\ell(\mu^i)=l_i}\mathfrak{F}^\lambda_{\mu,\mu^1,...,\mu^k},$$
where $\mathfrak{F}^\lambda_{\mu,\mu^1,...,\mu^k}$ is defined in \Cref{mathfrakF}. 
\cref{Thm3} can be reformulated as follows:

\begin{thm}\label{thm marginal sums}
For every $k\geq1$, we have
$$\tau_b^{(k)}(t,\textbf{p},\textbf{q},\underline{u_1},\underline{u_2},...,\underline{u_k})=\sum_{n\geq0}\sum_{\underset{l_1,...,l_k\geq1}{\lambda,\mu\vdash n}}\frac{p_\lambda q_\mu u_1^{l_1}...u_k^{l_k}}{z_\lambda(1+b)^{\ell(\lambda)}}\sum_{(\delta_0,...,\delta_{k-1})\in\mathfrak F^\lambda_{\mu,l_1,...,l_k}}b^{\vartheta_\lambda(\delta_0,...,\delta_{k-1})}.$$
 where $\vartheta_\lambda(\delta_0,...,\delta_{k-1})$ is a non-negative integer which is zero if and only if each one of the matchings $\delta_0,...,\delta_{k-1}$ is bipartite.
\end{thm}
The purpose of this section is to use the $b$-weight of rooted-constellations given in \cref{CD} in order to define a statistic $\vartheta$ on $k$-tuples of matchings that satisfies \cref{thm marginal sums}. We recall that in \cref{prop1} we have established a bijection between $k+2$-tuples of matchings and labelled $k$-constellations. The difficulty here is that the sums run over $k$-tuples of matchings (we recall that in definition of $\mathfrak{F}^\lambda_{\mu^0,...,\mu^k}$ we fix the matchings $\delta_{-1}$ to be  $\varepsilon$ and the matching $\delta_k$ to be $\delta_\lambda$; see \Cref{mathfrakF}). It turns out that the convenient objects to consider are the \textit{face-labelled constellations}. The purpose of Sections \ref{sec face-labelled}, \ref{sec bweights face-labelled}, and \ref{sec equivalence} is to introduce face-labelled constellations and define $b$-weights on them. In \cref{sec face-labelled and matchings} we will establish a bijection between $\mathfrak{F}^\lambda_{\mu^0,...,\mu^k}$ and face-labelled constellations.

\subsection{Face-labelled constellations.}\label{sec face-labelled}
Face-labelled maps were introduced in \cite{B21} in the case of bipartite maps, we give here an analog definition for constellations.
We say that a $k$-constellation $\mathbf{M}$ is \textit{face-labelled} if each face is rooted (with a marked oriented corner of color 0 or equivalently with a marked right-path), and for every $j>0$, the faces of degree $j$ are labelled \textit{i.e.} if $\mathbf{M}$ contains $m_j>0$ faces of degree $j$, these faces are labelled by $\{1,2,..,m_j\}$. Face-labelled constellations will be denoted with a hat: $\hat{\mathbf{M}}$. In each face, the marked corner or right-path is called the \textit{face-root}.
We say that a connected face-labelled $k$-constellation is \textit{oriented} if the underlying rooted constellation is orientable, and the orientations given by the face roots are consistent, see \cref{consistent Orientations}. Finally, we say that a connected face-labelled $k$-constellation $\hat{\mathbf{M}}$ is \textit{rooted} if the underlying constellation has a root $c$, such that the orientation of the root face (given by the definition of a face-labelled constellation above) is the same as the orientation induced by the root $c$. This constellation will be denoted $(\hat{\mathbf{M}},c)$.

\subsection{\texorpdfstring{$b$}{}-weights for connected face-labelled constellations.}\label{sec bweights face-labelled}

Once and for all, and for every connected rooted $k$-constellation $(\mathbf{M},c)$, we fix an orientation $O_{(\mathbf{M},c)}$ of the faces of $\mathbf{M}$ that satisfies the two following properties (see \cite[Section 5.1]{D17}):
\begin{itemize}
    \item The orientation of the root face is given by the root $c$.
    \item If $\mathbf{M}$ is orientable, then $O_{(\mathbf{M},c)}$ is the canonical orientation of the constellation, see \Cref{ssec Maps}.
\end{itemize}

\begin{defi}\label{MONFLR}
Let $(\hat{\mathbf{M}},c)$ be a connected rooted face-labelled constellation, and let $(\mathbf{M},c)$ be the underlying rooted constellation. We define the $b$-weight $\vartheta$ of $(\hat{\mathbf{M}},c)$ by $$\vartheta(\hat{\mathbf{M}},c):=\nu(\tilde{ \mathbf{M}},\tilde c)+r,$$  where $(\tilde{\mathbf{ M}},\tilde c)$ is the dual constellation of $(\mathbf{M},c)$ as defined in \cref{def duality}, $\nu(\tilde{\mathbf{M}},\tilde c)$ is the non-negative integer of \cref{CD},  and $r$ is the number of faces of $\hat{\mathbf{M}}$ whose orientation is different from the orientation given by $O_{(\mathbf{M},c)}$.
\end{defi}

\begin{rmq}\label{rmq oriented}
We note that $\vartheta(\hat{\mathbf{M}},c)=0$ if and only if $\hat{\mathbf{M}}$ is oriented. Moreover, for every connected rooted constellation $\mathbf{M}$ with face-type $\lambda$, we have
\begin{equation}\label{eqet}
\sum_{(\hat{\mathbf{M}},c)}b^{\vartheta(\hat{\mathbf{M}},c)}=z_\lambda(1+b)^{\ell(\lambda)-1}b^{\nu(\mathbf{\tilde M},\tilde c)},
\end{equation} 
where the sum is taken over all possible face-labellings of $(\mathbf{M},c)$. Indeed, we have $z_\lambda$ choices to label the faces of $(\mathbf{M},c)$ which have the same size and choose a corner of color 0 (which is not yet oriented) in each face. Besides, for each face other than the root face, we have to chose an orientation for the root corner (the orientation in the root face being fixed by the root $c$). The orientation consistent with $O_{(\mathbf{M},c)}$ contributes with 1 to the $b$-weight and the other orientation contributes with $b$, which gives us $1+b$ for each face different from the root face.
\end{rmq}

We now define $b$-weights for unrooted connected face-labelled constellations. These $b$-weights are given by different ways to root a face-labelled constellation.

\begin{defi}\label{MONFLU}
Let $\lambda$ be a partition of $n$ and let $\hat{\mathbf{M}}$ be a connected $k$-constellation of face-type $\lambda$. We define three $b$-weights on $\hat{\mathbf{M}}$:
\begin{enumerate}
    \item We root $\hat{\mathbf{M}}$ with $c_0$, the root of the face of maximal degree and which is labelled by 1. We define :
    $$\Vec{\rho}(\hat{\mathbf {M}}):=b^{\vartheta(\hat{\mathbf{M}},c_0)}.$$
    \item We take the average over all possible roots $c$ that lies in a face of maximal degree (we recall that the orientation of this root should be consistent with the orientation given by the face-root): $$\hat{\rho}(\hat{\mathbf{M}}):=\frac{1}{m\lambda_1}\sum_{c,\deg(f_c)=\lambda_1}b^{\vartheta(\hat{\mathbf{M}},c)},$$
    where $m:=m_{\lambda_1}(\lambda)$ is the number of faces of maximal degree.
    \item We take the average over all possible roots $c$ :
    $$\tilde{\rho}(\hat{\mathbf{M}}):=\frac{1}{n}\sum_{c}b^{\vartheta(\hat{\mathbf{M}},c)}.$$
\end{enumerate}
\end{defi}
\noindent Note that the $b$-weight $\vec{\rho}$ has the advantage of being integral, however it is a priory less symmetric than $\tilde{\rho}$. The purpose of the next subsection is to show that the $b$-weights $\vec \rho$, $\hat{\rho}$ and $\tilde{\rho}$ are equivalent when we sum over connected face-labelled $k$-constellations of a given marginal profile $(\lambda,\mu,l_1,...,l_k)$. 

\subsection{Equivalence between the three \texorpdfstring{$b$}{}-weights.}\label{sec equivalence}

We start by the equivalence between $\vec{\rho}$ and~$\hat{\rho}$.
\begin{lem}\label{vec-hat}
For every $k,n\geq1$ and
and $\lambda,\mu^0,...,\mu^k\vdash n$ we have
\begin{equation}\label{eq vec-hat}
    \sum_{\hat{\mathbf{M}}}\vec{\rho}(\hat{\mathbf{M}})=\sum_{\hat{\mathbf{M}}}\hat{\rho}(\hat{\mathbf{M}}),
\end{equation}
where the sums are taken over connected face-labelled $k$-constellation with profile $(\lambda,\mu^0,...,\mu^k)$.
\begin{proof}
We denote $m:=m_{\lambda_1}(\lambda)$, the number of parts in $\lambda$ of maximal size.
From \cref{MONFLR} and \cref{MONFLU}, we know that  $\vec{\rho}(\hat{\mathbf{M}})$ is of the form $b^rb^{\nu(\tilde{\mathbf {M}},\tilde c)}$. 
We rewrite the left-hand side of \cref{eq vec-hat} by putting together the terms having the same underlying rooted connected constellation $(\mathbf{M},c)$. With the same argument as in the proof of \cref{eqet}, for every rooted constellation $(\mathbf{M},c)$ with a root $c$ in a face of maximal degree we have
$$\sum_{\hat{\mathbf{M}}}\vec{\rho}(\hat{\mathbf{M}})=(1+b)^{\ell(\lambda)-1}\frac{z_\lambda}{m\lambda_1}b^{\nu(\tilde{\mathbf{M}},\tilde c)},$$
where the sum is taken over face-labelled constellations that can be obtained from $(\mathbf{M},c)$ by labelling its faces, with the condition that the root face is always labelled by 1 and rooted by $c$ (see \cref{MONFLU} item 1).
We deduce that the left-hand side of \cref{eq vec-hat} side is equal to 
$$(1+b)^{\ell(\lambda)-1}\frac{z_\lambda}{m\lambda_1}\sum_{(\mathbf{M},c)}b^{\nu(\tilde{\mathbf{M}},\tilde c)}$$
where the sum is taken over rooted connected $k$-constellations with profile $(\lambda,\mu^0,...,\mu^k)$, such that the root face has maximal degree $\lambda_1$.

\noindent On the other hand, we can rewrite the right-hand side of \cref{eq vec-hat} (using \cref{MONFLU} item 2) as follows
$$\sum_{(\hat{\mathbf{M}},c)}\frac{1}{m\lambda_1}b^{\vartheta(\hat{\mathbf{M}},c)},$$
where the sum is taken over face-labelled rooted constellation, for which the root is in a face of maximal degree. 

\noindent We use \cref{eqet} to conclude.
\end{proof}
\end{lem}

The link between the two $b$-weights $\hat{\rho}$ and $\tilde{\rho}$ is less obvious. We need a property of symmetry of the $b$-weight defined in \cite{CD20} on rooted connected constellations. We start by defining for every $s\geq1$ the series 
$$U_s:=(1+b)s\frac{\partial}{\partial q_s}\log(\tau^{(k)}_b)\hspace{0,4cm},\hspace{0,4cm}V_s:=(1+b)s\frac{\partial}{\partial p_s}\log(\tau^{(k)}_b).$$
We also define the operator $\pi$ that switches the variables $\mathbf{p}\leftrightarrow\mathbf{q}$ and $u_i\leftrightarrow u_{k+1-i}$ for $1\leq i\leq k$.
Since $\pi\tau^{(k)}_b=\tau^{(k)}_b$, we get $\pi U_s=V_s$. On the other hand, one has (see \cite[Corollary 5.9]{CD20})
\begin{equation}\label{eq U}
    U_s=q_s^{-1}\sum_{\underset{\deg(v_c)=s} {(\mathbf{M},c)} }t^{|\mathbf{M}|}\kappa(\mathbf{M})b^{\nu(\mathbf{M},c)},
\end{equation}
where the sum is taken over rooted connected $k$-constellation whose root vertex has degree $s$. Moreover, it is straightforward from \cref{rmq duality} that for every $k$-constellation $\mathbf{M}$ we have 
\begin{equation}\label{eq pi}
    \pi\big(\kappa(\mathbf{M})\big)=\kappa(\tilde{\mathbf{M}}),
\end{equation}
where $\tilde{\mathbf{M}}$ denotes the dual constellation  of $\mathbf{M}$.
Applying $\pi$ to \cref{eq U}, we get
\begin{equation}\label{eq Vs}
    V_s=p_s^{-1}\sum_{\underset{\deg(f_c)=s}{(\mathbf{M},c)}}t^{|\mathbf{M}|}\kappa(\mathbf{M})b^{\nu(\tilde{\mathbf{M}},\tilde c)}.
\end{equation}
We deduce the following lemma.
\begin{lem}\label{Sym}
Let $\lambda,\mu\vdash n$ and $l_1,...l_k\geq1$, and let $s\geq1$ such that $m:=m_s(\lambda)\geq1$. Then 
$$\frac{1}{n}\sum_{(\mathbf{M},c)}b^{\nu(\tilde{\mathbf{M}},\tilde c)}=\frac{1}{ms}\sum_{\underset{\deg(f_c)=s}{(\mathbf{M},c)}}b^{\nu(\tilde{\mathbf{M}},\tilde c)},$$
where the sums are taken over connected rooted $k$-constellations of marginal profile $(\lambda,\mu,l_1,...,l_k)$, with the condition that the root face has degree $s$ in the sum of the right-hand side.
\begin{proof}
From \cref{CD} we have 
$$(1+b)\log(\tau^{(k)}_b)=\sum_{(\textbf{M},c)}\frac{t^{|\textbf{M}|}}{|\textbf{M}|} \kappa(\textbf{M}) b^{\nu(\mathbf{M}, c)}.$$
Applying $\pi$ on the last equality, we get
$$(1+b)\log(\tau^{(k)}_b)=\sum_{(\textbf{M},c)}\frac{t^{|\textbf{M}|}}{|\textbf{M}|} \kappa(\mathbf{M}) b^{\nu(\mathbf{\tilde M},\tilde c)}.$$
We deduce then that the coefficient of the monomial $t^np_\lambda q_{\mu^0}u_1^{l_1}...u_k^{l_k}$ in $p_sV_s$, is given by
$$\frac{ms}{n}\sum_{(\mathbf{M},c)}b^{\nu(\mathbf{\tilde M},\tilde c)},$$
where the sum is taken over connected rooted $k$-constellations of marginal profile $(\lambda,\mu,l_1,...,l_k)$.
On the other hand, using \cref{eq Vs} we get that this coefficient is also equal to 
$$\sum_{\underset{\deg(f_c)=s}{(\mathbf{M},c)}}b^{\nu(\mathbf{\tilde M},\tilde c)},$$
where the sum is taken over connected rooted $k$-constellations of marginal profile $(\lambda,\mu,l_1,...,l_k)$ with the condition that the root face has degree $s$,
which finishes the proof.
\end{proof}
\end{lem}
This lemma has the following interpretation : conditioning the root to be in a face of a given degree does not affect the $b$-weight obtained when summing over constellations of a given marginal profile.
We deduce the following corollary that establishes the equivalence claimed between $\hat{\rho}$ and $\tilde{\rho}$:
\begin{cor}\label{hat-tilde}
Let $\lambda,\mu\vdash n$, and $l_1,...,l_k\geq1$. Then we have 
$$\sum_{\hat{\mathbf{M}}}\hat\rho(\hat{\mathbf{M}})=\sum_{\hat{\mathbf{
M}}}\Tilde{\rho}(\hat{\mathbf{M}}),$$
where the sums run over connected face-labelled $k$-constellation of marginal profile $(\lambda,\mu,l_1,...,l_k)$.
\begin{proof}
We apply \cref{Sym} for $s=\lambda_1$ and we multiply the equation  by $z_\lambda(1+b)^{\ell(\lambda)-1}$. Using \cref{eqet} we obtain:
$$\frac{1}{n}\sum_{(\hat{\mathbf M},c)}b^{\vartheta(\hat{\mathbf{M}},c)}=\frac{1}{m\lambda_1}\sum_{\underset{\deg(f_c)=\lambda_1}{(\hat{\mathbf{M}},c)}}b^{\vartheta(\hat{\mathbf{M}},c)},$$
where $m:=m_{\lambda_1}(\lambda)$, which finishes the proof.
\end{proof}
\end{cor}

\subsection{Extension to disconnected face-labelled constellations.}\label{ssec disconnected}
We extend multiplicatively the $b$-weight $\vec{\rho}$ to disconnected constellations. More precisely, if $\hat{\mathbf{M}}$ is a disconnected face-labelled constellation and $\hat{\mathbf{M}}_i$ is a connected component of $\hat{\mathbf{M}}$, then it can be considered as a face-labelled constellation where the labelling of the faces having the same degree in $\hat{\mathbf{M}}_i$ is inherited from their labelling in $\hat{\mathbf{M}}$.  This allow us to define $\vec \rho(\hat{\mathbf{M}})$ as the product over all its connected components of $\vec \rho(\hat{\mathbf{M}}_i)$, where  $\vec \rho(\hat{\mathbf{M}}_i)$ is given by \cref{MONFLU} item 1. 
\begin{rmq}\label{rmq disconnected}
By definition $\hat{\mathbf{M}}$ is oriented if and only if each one of its connected components is oriented. Hence, we can deduce from \cref{rmq oriented} that $\vec{\rho}(\mathbf{\hat{M}})$ is a monomial, and it equals 1 if and only if $\mathbf{\hat{M}}$ is oriented. Hence $\vec{\rho}$ is an integral $b$-weight on face-labelled constellations.
\end{rmq}
The following lemma establishes the connection between the generating functions of connected and disconnected constellations. It is a variant of the exponential formula in the combinatorial class theory. However, one has to take care of the multiplicities since we have a separate labelling for each size of faces. We give here the proof in completeness. 

\begin{lem}\label{disconnected}
For every $k \geq 1$, we have
\begin{multline*}
    \sum_{n\geq1}t^n\sum_{\lambda,\mu\vdash n}\sum_{l_1,..,l_k\geq1}\frac{p_\lambda q_\mu u_1^{l_1}...u_k^{l_k}}{z_\lambda(1+b)^{\ell(\lambda)}}\sum_{\hat{\mathbf{M}}}\vec \rho(\hat{\mathbf{M}})\\
    =\exp\Big(\sum_{n\geq1}t^n\sum_{\lambda,\mu\vdash n}\sum_{l_1,..,l_k\geq1}\frac{p_\lambda q_\mu u_1^{l_1}...u_k^{l_k}}{z_\lambda(1+b)^{\ell(\lambda)}}\sum_{\hat{\mathbf{M}}\text{ connected}}\vec \rho(\hat{\mathbf{M}})\Big),
\end{multline*}
where the last sum is taken each time on face-labelled constellations of marginal profile $(\lambda,\mu,l_1,...,l_k)$.
\begin{proof}
When we develop the exponential of the right-hand side, we obtain a sum over tuples of connected face-labelled constellations. Let $\mathbf{\hat{M}}_1$,...,$\mathbf{\hat{M}}_r$ be a list of $r$ connected face-labelled constellations, with face-types $\lambda^{(1)}$,...,$\lambda^{(r)}$. We define $\lambda:=\bigcup\limits_{i=1}^r\lambda_i$. Taking the disjoint union of the constellations $\mathbf{\hat{M}}_1$,...,$\mathbf{\hat{M}}_r$, we obtain a constellation of face-type $\lambda$. 
In such operations, we deal with the labellings as in the theory of labelled
combinatorial classes \cite[Chapter II]{FS09}; namely for every $j$ such that $m_j(\lambda)>0$,  we consider all the ways to relabel the faces of degree $j$ of $\mathbf{\hat{M}}_1$,...,$\mathbf{\hat{M}}_r$ in an increasing way such that their label sets become disjoint and the union of their
label sets is $\llbracket m_j(\lambda)\rrbracket$. So we have $\binom{m_j(\lambda)}{m_j(\lambda^1),...,m_j(\lambda^r)}$ choices to relabel the faces of degree $j$, and 
$$\frac{z_\lambda}{z_{\lambda^1}...z_{\lambda^r}}=\prod_{j}\binom{m_j}{m_j(\lambda^1),,...,m_j(\lambda^r)}$$ choices to relabel all the faces of $\bigcup\limits_{i=1}^r\mathbf{\hat{M}}_i$ to obtain a face-labelled constellation $\hat{\mathbf{M}}$. 
Finally, we notice that the marking and the quantity $(1+b)^{\ell(\lambda)}$ are multiplicative which concludes the proof.
\end{proof}
\end{lem}

\subsection{Face-labelled constellations and Matchings }\label{sec face-labelled and matchings}
Let $\hat{\mathbf{M}}$ be a face-labelled constellation of face-type $\lambda$. We describe a canonical way to obtain a labelled constellation $\Check{\mathbf{M}}$ from $\hat{\mathbf{M}}$ that will be useful in the next proposition. We start by defining the following order on $\mathcal{A}_n$: $\hat{1}<1<2...<\hat{n}<n$. We label the right-paths starting from faces of highest degree and smallest label: We start from the face of degree $\lambda_1$ and label 1. We travel along this face starting from the right-path preceding the root corner, and we attribute to each right-path the smallest label not yet used. We restart with the face of highest degree and smallest label whose right-paths are not yet labelled. This will be called the \textit{canonical labelling} of the face-labelled constellation $\hat{\mathbf{M}}$. Note that for every face-labelled constellation $\hat{\mathbf{M}}$ of face-type $\lambda$, the matchings $\delta_{-1}$ and $\delta_k$ associated to the canonical labelling $\Check{\mathbf{M}}$ (as in \cref{Matchings}) satisfy $\delta_{-1}=\varepsilon$ and $\delta_k=\delta_\lambda$, where $\delta_\lambda$ is the matching defined in \cref{def delta lambda}.
\medskip

\noindent We now prove the following proposition that establishes a correspondence between face-labelled $k$-constellations and $k$-tuples of matchings.

\begin{thm}\label{prop2}
For $\lambda,\mu^0,...,\mu^k\vdash n$, there exists a bijection $\mathcal{\hat M}$ between
face-labelled $k$-constellations with profile $(\lambda,\mu^0,...,\mu^k)$ and $\mathfrak{F}^\lambda_{\mu^0,...,\mu^k}$. Moreover, a face-labelled constellation $\mathbf{\hat{M}}$ is oriented if and only if $\mathcal{\hat M}(\mathbf{\hat{M}})$ is a $k$-tuple of bipartite matchings.
\begin{proof}
Let $\mathbf{\hat{M}}$ be a face-labelled constellation, and let $\mathbf{\Check{M}}$ be its canonical labelling.
By construction of the canonical labelling, the matchings associated to $\mathbf{\Check{M}}$ by the bijection of \cref{prop1} is of the form $\mathcal{M}(\mathbf{\Check{M}})=(\varepsilon,\delta_0,...,\delta_{k-1},\delta_\lambda)$. Then we define 
$\mathcal{\hat M}(\mathbf{\hat{M}}):=(\delta_0,...,\delta_{k-1})$.

\noindent Conversely, let $(\delta_0,...,\delta_{k-1})\in\mathfrak{F}^\lambda_{\mu^0,...,\mu^k}$. The bijection of \cref{prop1} gives us a labelled $k$-constellation $\mathbf{\Check{M}}:=\mathcal{M}^{-1}(\varepsilon,\delta_0,...,\delta_{k-1},\delta_\lambda)$.
Then $\mathcal{\hat M}^{-1}(\delta_0,...,\delta_{k-1})$ is the face-labelled constellation having $\mathbf{\Check{M}}$ as a canonical labelling.

We now prove the second part of the proposition. Let $\mathbf{\hat{M}}$ be a face-labelled constellation. We start by the following remark: When we travel along the boundary of each face of $\mathbf{\hat{M}}$ in the orientation induced by its root, a right-path $H$ in $\mathbf{\hat{M}}$ is traversed from the corner of color 0 to the corner of color $k$ (respectively from the corner of color $k$ to the corner of color 0) if it has a label in the first class of $\mathcal{A}_n$ (respectively in second class) with respect to the canonical labelling of $\mathbf{\hat{M}}$.
Indeed, this property is clear for the root right-path of each face. Moreover, this can be extended to the other right-paths since when we travel along a face we alternate right-paths traversed form 0 to $k$ and right-paths traversed from $k$ to 0, and when we traverse a connected component of $G(\varepsilon,\delta_\lambda)$ we alternate labels of first and second class.

\noindent We recall that $\mathbf{\hat{M}}$ is oriented if and only if the orientations induced by the faces roots are consistent as in \cref{consistent Orientations}. Note that the orientations of the faces of $\mathbf{\hat{M}}$ are consistent from either side of edges of color $(i,i+1)$ if and only if each two right-paths having an edge of color $(i,i+1)$ in common are traversed in opposite ways. By the previous remark, this is equivalent to say that $\delta_i$ is bipartite. In particular $\hat{\mathbf{M}}$ is oriented if and only if each one of the matchings $\delta_0$,...,$\delta_{k-1}$ is bipartite.
\end{proof}
\end{thm}

\subsection{A statistic \texorpdfstring{$\vartheta$}{} for elements of \texorpdfstring{$\mathfrak F^\lambda_{\mu,l_1,...,l_k}$}{} and proof of \texorpdfstring{\cref{Thm3}}. }
\begin{defi}\label{def-lem}
Let $\lambda,\mu\vdash n$ and $l_1,...,l_k\geq1$. For each $(\delta_1,...,\delta_k)\in\mathfrak F^\lambda_{\mu,l_1,...,l_k}$, we define the non-negative integer
$\vartheta_\lambda(\delta_0,..,\delta_{k-1})$ such that
$$\vec\rho(\hat{\mathbf{M}})=b^{\vartheta_\lambda(\delta_0,..,\delta_{k-1})},$$
where $\hat{\mathbf{M}}$ is the face-labelled constellation associated to $(\delta_1,...,\delta_k)$ by the bijection of \cref{prop2}. 
\end{defi}
Since the bijection of \cref{prop2} ensures that  $\hat{\mathbf{M}}$ is oriented if and only if the matchings $\delta_0,...,\delta_{k-1}$ are bipartite, we note that $\vartheta_\lambda(\delta_0,..,\delta_{k-1})$  is equal to zero if and only if each one of the matchings $\delta_1$,...,$\delta_k$ is bipartite.

\begin{proof}[Proof of \cref{thm marginal sums}]
From \cref{CD} we have

$$\Psi_b^{(k)}(t,\mathbf{p},\mathbf{q},\underline{u_1},\underline{u_2},...,\underline{u_k})=t\frac{\partial}{\partial t}\sum_{n\geq1}\frac{t^n}{n}\sum_{(\mathbf{M},c)} \kappa(\mathbf{M}) b^{\nu(\mathbf{M},c)},$$
where the second sum runs over connected rooted constellations of size $n$.
Applying the operator $\pi$ on the last equation and using \Cref{eq pi}  we get
$$\Psi_b^{(k)}(t,\mathbf{p},\mathbf{q},\underline{u_1},\underline{u_2},...,\underline{u_k})=t\frac{\partial}{\partial t}\sum_{n\geq1}\frac{t^n}{n}\sum_{(\mathbf{M},c)} \kappa(\mathbf{M}) b^{\nu(\tilde{\mathbf{M}},c)}.$$
 Using \cref{eqet} and \Cref{MONFLU} item \textit{(3)}, we obtain 
\begin{align*}
    \Psi_b^{(k)}(t,\mathbf{p},\mathbf{q},\underline{u_1},\underline{u_2},...,\underline{u_k})
    &=t\frac{\partial}{\partial t}\sum_{n\geq1}\frac{t^n}{n}\sum_{\lambda,\mu\vdash n}\sum_{l_1,..,l_k\geq1}\frac{p_\lambda q_\mu u_1^{l_1}...u_k^{l_k}}{z_\lambda(1+b)^{\ell(\lambda)-1}}\sum_{(\hat{\mathbf{M}},c)}b^{\vartheta(\hat{\mathbf{M}},c)}\\
    &=t\frac{\partial}{\partial t}\sum_{n\geq1}t^n\sum_{\lambda,\mu\vdash n}\sum_{l_1,..,l_k\geq1}\frac{p_\lambda q_\mu u_1^{l_1}...u_k^{l_k}}{z_\lambda(1+b)^{\ell(\lambda)-1}}\sum_{\hat{\mathbf{M}}}\tilde \rho(\hat{\mathbf{M}}),
\end{align*}

\noindent where the last sums are taken over connected constellations of marginal profile $(\lambda,\mu,l_1,...,l_k)$, rooted in the first equation and unrooted in the second one. This can be rewritten using \cref{vec-hat} and \cref{hat-tilde} as follows 
$$\Psi_b^{(k)}(t,\mathbf{p},\mathbf{q},\underline{u_1},\underline{u_2},...,\underline{u_k})
=t\frac{\partial}{\partial t}(1+b)\sum_{n\geq1}t^n\sum_{\lambda,\mu\vdash n}\sum_{l_1,..,l_k\geq1}\frac{p_\lambda q_\mu u_1^{l_1}...u_k^{l_k}}{z_\lambda(1+b)^{\ell(\lambda)}}\sum_{\hat{\mathbf{M}}\text{ connected}}\vec \rho(\hat{\mathbf{M}}).$$

Applying \cref{disconnected}, we obtain 
$$\Psi_b^{(k)}(t,\mathbf{p},\mathbf{q},\underline{u_1},\underline{u_2},...,\underline{u_k})=t\frac{\partial}{\partial t}(1+b)\log\Big(\sum_{n\geq1}t^n\sum_{\lambda,\mu\vdash n}\sum_{l_1,..,l_k\geq1}\frac{p_\lambda q_\mu u_1^{l_1}...u_k^{l_k}}{z_\lambda(1+b)^{\ell(\lambda)}}\sum_{\hat{\mathbf{M}}}\vec \rho(\hat{\mathbf{M}})\Big),$$
where the last sum runs over constellations of marginal profile $(\lambda,\mu,l_1,...l_k)$, connected or not. Comparing the last equation with \cref{eqPsi}, we deduce that 
$$\tau_b^{(k)}(t,\mathbf{p},\mathbf{q},\underline{u_1},\underline{u_2},...,\underline{u_k})= \sum_{n\geq1}t^n\sum_{\lambda,\mu\vdash n}\sum_{l_1,..,l_k\geq1}\frac{p_\lambda q_\mu u_1^{l_1}...u_k^{l_k}}{z_\lambda(1+b)^{\ell(\lambda)}}\sum_{\hat{\mathbf{M}}}\vec \rho(\hat{\mathbf{M}}).$$
Using the bijection of \cref{prop2} and \cref{def-lem}, the last equation can be rewritten as follows:
$$\tau_b^{(k)}(t,\textbf{p},\textbf{q},\underline{u_1},\underline{u_2},...,\underline{u_k})=\sum_{n\geq0}\sum_{\underset{l_1,...,l_k\geq1}{\lambda,\mu\vdash n}}\frac{p_\lambda q_\mu u_1^{l_1}...u_k^{l_k}}{z_\lambda(1+b)^{\ell(\lambda)}}\sum_{(\delta_0,...,\delta_{k-1})\in\mathfrak F^\lambda_{\mu,l_1,...,l_k}}b^{\vartheta_\lambda(\delta_0,...,\delta_{k-1})}.$$
Since $\vartheta_\lambda$ has the properties required in \cref{thm marginal sums}, this concludes the proof.
\end{proof}

\section{Application: Lassale's conjecture for rectangular partitions}\label{sec Lassale conjecture}

In this section, we consider the function $\tau_b^{(1)}$ with the following specialisations; $\mathbf{q}^{(0)}=\underline{-r\alpha}:=(-r\alpha,-r\alpha,...)$ and $\mathbf{q}^{(1)}=\underline{q}:=(q,q,...)$. We also replace $t$ by $-t$.
We recall that for $k=1$, $k$-constellations correspond to bipartite maps; these are maps with vertices colored in white and black, such that each edge separates two vertices of different colors.  For every face-labelled bipartite map $\mathbf{M}$, we denote by $\w^{(\alpha)}(\hat{\mathbf{M}},q,r)$ the quantity:
$$\w^{(\alpha)}(\hat{\mathbf{M}},q,r):=(-1)^{|\mathbf{M}|}(-r\alpha)^{|\mathcal{V}_\circ(\hat{\mathbf{M}})|}q^{|\mathcal{V}_\bullet(\hat{\mathbf{M}})|},$$
where $|\mathcal{V}_\circ(\hat{\mathbf{M}})|$ (respectively $|\mathcal{V}_\bullet(\hat{\mathbf{M}})|$) denotes the number of white (respectively black) vertices of $\hat{\mathbf{M}}$.
The generating series of bipartite maps can be written as follows (see \Cref{thm marginal sums}):
\begin{equation}\label{gen series bip maps}
  \tau^{(1)}_b(-t,\mathbf{p},\underline{q},\underline{-r\alpha})=\sum_{\mu\in\mathcal{P}}\frac{p_\mu t^{|\mu|}}{z_\mu \alpha^{\ell(\mu)}}\sum_{\hat{\mathbf{M}}}\Vec{\rho}(\hat{\mathbf{M}})\w^{(\alpha)}(\hat{\mathbf{M}},q,r),  
\end{equation}
where the second sum runs over face-labelled bipartite maps of face-type $\mu$.

The key step of the proof is the following lemma, that gives an expression of Jack polynomials associated to partitions with rectangular shape in terms of the function $\tau^{(1)}_b$.
\begin{lem}\label{lem recJack}
For every partition $\lambda=(q\times r)\vdash n$, we have 
\begin{equation}\label{eq recJack}
  J_\lambda^{(\alpha)}=[t^n]\tau^{(1)}_b(-t,\mathbf{p},\underline{q},\underline{-r\alpha}),
\end{equation}
where $[.]$ denotes the extraction symbol with respect to the variable $t$.
\begin{proof}
Recall that 
\begin{equation}\label{eq rectau}
  [t^n]\tau^{(1)}_b(-t,\mathbf{p},\underline{q},\underline{-r\alpha})=(-1)^n\sum_{\nu\vdash n}\frac{J^{(\alpha)}_\nu(\mathbf{p})J^{(\alpha)}_\nu(\underline{q}) J_\nu^{(\alpha)}(\underline{-\alpha r})}{j_\nu^{(\alpha)}}  
\end{equation}
Let $\Box_0$ be fixed box, and let $u:=-c_\alpha(\Box_0)$ be the opposite of its $\alpha$-content, see \Cref{subsec Partitions}. Using \Cref{Jack formula}, we can see that
$J_\nu(\underline{u})=0$ if and only if $\Box_0\in \nu$. In particular, the partitions  that contribute in the sum of \Cref{eq rectau} are the partitions that do not contain the boxes of coordinates $(r+1,1)$ and $(1,q+1)$. The only partition of size $n$ that fulfills this condition is the partition $\lambda$. Hence 
$$[t^n]\tau^{(1)}_b(-t,\mathbf{p},\underline{q},\underline{-r\alpha})=(-1)^n\frac{J^{(\alpha)}_\lambda(\mathbf{p})J^{(\alpha)}_\lambda(\underline{q}) J_\lambda^{(\alpha)}(\underline{-\alpha r})}{j_\lambda^{(\alpha)}}.$$
Moreover, one can check that (see \Cref{Jack formula} and \Cref{eq j alpha})
$$(-1)^n\frac{J^{(\alpha)}_\lambda(\underline{q}) J_\lambda^{(\alpha)}(\underline{-\alpha r})}{j_\lambda^{(\alpha)}}=1,$$
which concludes the proof.
\end{proof}
\end{lem}

Note that the previous lemma can be generalized as follows; let $1=q_1<q_2<...q_s$ and $r_1>,...>r_s=1$ be two sequences of positive integers. We define the boxes $\Box_i:=(r_i,q_i)$ for $1\leq i\leq s$, and let $u_i:=-c_\alpha(\Box_i)$ be the respective opposite of their $\alpha$-contents. We define $\lambda$ as the partition of maximal size that does not contain any one of the boxes $\Box_i$, and let $n$ be the size of $\lambda$. One can write 
\begin{equation}
  J_\lambda=\frac{[t^n]\tau^{(k-1)}_b\left(t,\mathbf{p},\underline{u_1}...,\underline{u_k}\right)}{[t^n]\tau^{(k-2)}_b\left(t,\underline{u_1}...,\underline{u_k}\right)}.  
\end{equation}
However, we were not able to use the previous equation for partitions $\lambda$ which are not rectangular to get polynomiality information on coefficients $\theta^{(\alpha)}_\mu(\lambda)$.

The purpose of the following lemma is to explain how to add faces of degree 1 on $b$-weighted bipartite maps. We will need a variant of \Cref{gen series bip maps} where we replace $\Vec{\rho}$ by another $b$-weight on face-labelled maps $\Vec{\rho}_{SYM}$  that we now define. As noticed in \Cref{rmq thmCD20}, the $b$-weight $b^{\nu(\mathbf{M},c)}$ that we consider in \Cref{sec3} is not the only one that satisfies \Cref{CD}. We consider now the $b$-weight $\rho_{SYM}(\mathbf{M},c)$ defined in \cite[Remark 3]{CD20}, which is not integral but has more symmetry properties that will be useful in the proof of \Cref{lem add 1-face}. We define $\Vec{\rho}_{SYM}$ as the $b$-weight on face-labelled bipartite maps obtained in \Cref{ssec disconnected}  when we replace $b^{\nu(\mathbf{M},c)}$ by $\rho_{SYM}(\mathbf{M},c)$ in \Cref{sec3} (see also \Cref{MONFLR} and \Cref{MONFLU} item(1)). With the same arguments used in \Cref{sec3}, one can check that \Cref{gen series bip maps} also  holds for $\Vec{\rho}_{SYM}$.

\begin{lem}\label{lem add 1-face}
For every partition $\mu\vdash m$, such that $m_1(\mu)=0$, and $\lambda=(q\times r)\vdash n\geq m$, we have
\begin{equation}\label{eq add 1}
[p_{\mu\cup1^{n-m}}t^n]\tau^{(1)}_b(-t,\mathbf{p},\underline{q},\underline{-r\alpha})=[p_\mu t^m]\tau^{(1)}_b(-t,\mathbf{p},\underline{q},\underline{-r\alpha}),
\end{equation}
where $\mu\cup1^\ell$ denotes the partition obtained by adding $\ell$ parts equal to 1 to the partition $\mu$.
\begin{proof}
We start by proving that for every partition $\nu\vdash \ell$ we have the following equation;
\begin{equation}\label{eq add 1-face}
  2(m_1(\nu)+1)[p_{ \nu\cup 1} t^{\ell+1}]\tau^{(1)}_b(-t,\mathbf{p},\underline{q},\underline{-r\alpha})=2(n-\ell)[p_{\nu} t^\ell]\tau^{(1)}_b(-t,\mathbf{p},\underline{q},\underline{-r\alpha}).  
\end{equation}
Using \Cref{gen series bip maps} for the $b$-weight $\Vec{\rho}_{SYM}$ introduced above, we can see that the two terms of the previous equation are generating series of bipartite maps. Hence, the last equality can be rewritten as follows;
\begin{equation}\label{eq add 1-face 2}
2\sum_{\hat{\mathbf{M}}}\Vec{\rho}_{SYM}(\hat{\mathbf{M}})\w^{(\alpha)}(\hat{\mathbf{M}})= 2\alpha(n-l)\sum_{\hat{\mathbf{M}}}\Vec{\rho}_{SYM}(\hat{\mathbf{M}})\w^{(\alpha)}(\hat{\mathbf{M}})    
\end{equation}
where the sums run over face-labelled bipartite maps, of face-type $\nu\cup 1$ in the left hand-side and $\nu$ in the right hand-side. The factor 2 in the left hand side of the last equation will be interpreted as marking an edge-side in the face of degree 1 with the highest label of each face-labelled bipartite map $\hat{\mathbf{M}}$ of face-type $\nu\cup 1$.
Such a map can be obtained by adding an edge $e$ with a marked side to a bipartite map $\hat{\mathbf{M}}$ of face-type $\nu$, such that the marked side is in a face of degree 1 in the map $\hat{\mathbf{M}}\cup\{e\}$. 

\noindent In the following, we show that this corresponds to the right-hand side of \Cref{eq add 1-face 2}. Let $\hat{\mathbf{M}}$ be a map of face-type $\nu$. We have two ways to add such an edge $e$ with a marked side to $\hat{\mathbf{M}}$:
\begin{itemize}
    \item We add an isolated edge with a marked side. We chose the highest label for the face of degree 1 that we form by adding $e$. We thus obtain a face-labelled map. In this case we have : $\w(\hat{\mathbf{M}}\cup\{e\})=2n\alpha.\w(\hat{\mathbf{M}})$; the black vertex has weight $q$, the white $-r\alpha$, and we multiply by $-1$ for adding an edge. Finally we have two choices for the marked edge-side. 
    \item We choose a side of an edge $s$ to which we add the marked side of the edge $e$ in order to form a face of degree 1. Since the map is of size $\ell$ we have $2\ell$ choices for the edge-side $s$. Once $s$ is fixed, we chose the highest label for the face of degree 1 formed by adding $e$. Since we have two choices of the orientation of this face, we obtain two face-labelled maps of face-type $\nu\cup 1$, that we denote $\hat{\mathbf{M}}_1$ and $\hat{\mathbf{M}}_2$. They satisfy $\w(\hat{\mathbf{M}}_1)=\w(\hat{\mathbf{M}}_2)=-\w(\hat{\mathbf{M}})$.
\end{itemize}
On the other hand, we claim that the $b$-weight $\Vec{\rho}_{SYM}$ defined above has the following property: if $e$ is an edge that we add to a bipartite map $\hat{\mathbf{M}}$ to form a face of degree then we have:
\begin{itemize}
    \item $\Vec{\rho}_{SYM}(\hat{\mathbf{M}}\cup\{e\})=\Vec{\rho}_{SYM}(\hat{\mathbf{M}})$, if $e$ is an isolated edge.
    \item $\Vec{\rho}_{SYM}(\hat{\mathbf{M}}_1)+\Vec{\rho}_{SYM}(\hat{\mathbf{M}}_2)=\alpha\Vec{\rho}_{SYM}(\hat{\mathbf{M}})$, if $e$ is added on an edge side of $\hat{\mathbf{M}}$, where $\hat{\mathbf{M}}_1$ and $\hat{\mathbf{M}}_2$ are as above.
\end{itemize}
Let us explain how to obtain this property. As explained above $\Vec{\rho}_{SYM}$ is obtained from $\rho_{SYM}$ by duality (see \Cref{MONFLR} and \Cref{MONFLU} item (1)). Notice that adding a face of degree 1 on a map is equivalent to adding a white vertex of degree 1 on the dual map. But such operation does not affect the $b$-weight $\rho_{SYM}$ (this is clear from the combinatorial model used in \cite{CD20} and the definition of $\rho_{SYM}$ \cite[Remark 3]{CD20}). Finally, observe that when $e$ is not an isolated edge, one of the possible orientations of the added face does not affect the $b$-weight $\rho_{SYM}(\hat{\mathbf{M}})$, and for the second one  $\rho_{SYM}(\hat{\mathbf{M}})$ is multiplied by $b$ (see \Cref{MONFLR}). This concludes the proof of the previous property and thus the proof of \Cref{eq add 1-face}.
 
\noindent Using \Cref{eq add 1-face}, we prove by induction on $\ell$ that 
$$[p_{\mu\cup 1^{\ell-m}} t^\ell]\tau^{(1)}_b(-t,\mathbf{p},\underline{q},\underline{-r\alpha})=\binom{n-m}{\ell-m}[p_\mu t^m]\tau^{(1)}_b(-t,\mathbf{p},\underline{q},\underline{-r\alpha}).$$
This gives \Cref{eq add 1} when $\ell=n$.
\end{proof}
\end{lem}
\begin{rmq}
Observe that \Cref{eq add 1-face} can be rewritten as follows  
$$\frac{\partial}{\partial p_1}\tau^{(1)}_b(-t,\mathbf{p},\underline{q},\underline{-r\alpha})=-\left[\frac{(-r\alpha)q}{\alpha}+\frac{t\partial}{\partial t}\right]\tau^{(1)}_b(-t,\mathbf{p},\underline{q},\underline{-r\alpha}).$$
This result can also be obtained using differential equation proved in \cite{CD20}. It is a $b$-deformation of the first Virasoro constraint known in the cases $b \in \{0,1\}$ (see \cite{KZ15} and \cite[Equation (17)]{CD20}).
\end{rmq}

\noindent We now prove the main result of this section.
\begin{thm}\label{Thm Lassale conjecture}
For every partition $\mu\vdash m\geq1$ such that $m_1(\mu)=0$, $(-1)^m z_\mu\theta^{(\alpha)}_\mu(q,r)$ is a polynomial in $(q,-r,b)$ with non-negative integer coefficients. More precisely, we have
\begin{equation}\label{eq thm theta}
  z_\mu\theta^{(\alpha)}_\mu(q,r)=\sum_{\hat{\mathbf{M}}}\Vec{\rho}(\hat{\mathbf{M}})\frac{\w^{(\alpha)}(\hat{\mathbf{M}},q,r)}{\alpha^{\ell(\mu)}},  
\end{equation}
where the sum is taken over face-labelled bipartite maps of face-type $\mu$.
\begin{proof}

We start by proving  \Cref{eq thm theta}. Since $m_1(\mu)=0$, we have that $\theta^{(\alpha)}_\mu(\lambda)=\theta^{(\alpha)}_{\mu\cup 1^{n-m}}(\lambda)$. Applying \Cref{lem recJack} and \Cref{lem add 1-face}, we get
$$\theta^{(\alpha)}_\mu(\lambda)=[p_{\mu\cup 1^{n-m}}t^n]\tau^{(1)}_b(-t,\mathbf{p},q,-r\alpha)=[p_{\mu}t^m]\tau^{(1)}_b(-t,\mathbf{p},q,-r\alpha).$$
Using \Cref{gen series bip maps}, this leads to \Cref{eq thm theta}.
\medskip

\noindent Let us now prove that \Cref{eq thm theta} implies the positivity and the integrality of the coefficients of $(-1)^m z_\mu\theta^{(\alpha)}_\mu(q,r)$. Since $\w^{(\alpha)}(\hat{\mathbf{M}},q,r)$ is a polynomial in $(q,-r,b)$ with non-negative integer coefficients, it suffices to eliminate the the term $\alpha^{\ell(\mu)}$ that appears in the denominator of the right-hand side of \Cref{eq thm theta}. We say that a bipartite map is \textit{weakly face-labelled} if it is obtained from a face-labelled bipartite map for which we keep the labelling of faces with same degree, but we forget the orientation of all the faces except for the face of maximal degree and smallest label in each connected component. For such a map, we have a natural notion of rooting for every connected component given by the root of the face of maximal degree and minimal label. Using a variant of \Cref{eqet}, \Cref{eq thm theta} can be rewritten as follows
\begin{equation}\label{sep alpha b}
  z_\mu\theta^{(\alpha)}_\mu(q,r)=\sum_{\mathbf{M}}\prod_{(\mathbf{M}_i,c_i)}b^{\nu(\tilde{\mathbf{M}}_i,c_i)}\frac{\w^{(\alpha)}(\mathbf{M}_i,q,r)}{\alpha},  
\end{equation}
where the sum is taken over weakly face-labelled bipartite maps $\mathbf{M}$ with face-type $\mu$ and the product runs over the connected components of $\mathbf{M}$ rooted as explained above. To conclude, notice that it is direct from the definition that $\w(\mathbf{M}_i,q,r)$ is divisible by $\alpha$. 
\end{proof}
\end{thm}

\begin{rmq}
As noticed by Lassale \cite[Conjecture 1, item \textit{(iii)}]{Las08}, $z_\mu$ is the good normalization to obtain integrality in \Cref{Thm Lassale conjecture}. Indeed, if $\mu\vdash m$
$$[q^m](-1)^mz_\mu\theta^{(\alpha)}_{\mu}(q,r)=(-r)^{\ell(\mu)},$$
where $[.]$ denotes the extraction symbol with respect to the variable $q$. To see this, observe that the only face-labelled bipartite map that contributes to the monomial $q^m$ in \Cref{eq thm theta} is the map of size $m$ and face-type $\mu$ that contains $m$ black vertices.
\begin{rmq}
Lassale suggested that the coefficients $\theta^{(\alpha)}_{\mu}(q,r)$ have a natural expression as a positive polynomial in both variables $\alpha$ and $b$ (see \cite[Conjecture 2]{Las08}). Such an expression can be obtained from \Cref{sep alpha b} by considering the two terms 
$\prod\limits_{(\mathbf{M}_i,c_i)}b^{\nu(\tilde{\mathbf{M}}_i,c_i)}$ and $\prod\limits_{(\mathbf{M}_i,c_i)}\frac{\w^{(\alpha)}(\mathbf{M}_i,q,r)}{\alpha}$. This expression in $\alpha$ and $b$ and the one given by Lassale in \cite{Las08} are related  but not the same.
\end{rmq}
\end{rmq}

\section{Generalized coefficients \texorpdfstring{$c^\lambda_{\mu^0,...,\mu^k}$}{}}\label{sec Top degree}

In this section, we state some properties of the coefficients $c^\lambda_{\mu^0,...,\mu^k}$, and we give a new proof for a combinatorial interpretation of the top degree part in these coefficients. This part is also related to the evaluation of $c^\lambda_{\mu^0,...,\mu^k}$ at $b=-1$, see \cref{cor b=-1}.

\subsection{General properties of \texorpdfstring{$c^\lambda_{\mu^0,...,\mu^k}$}{}}

We start by a multiplicativity property of the coefficients $c^\lambda_{\mu^0,...,\mu^k}$ due to Chapuy and Do\l{}{\k{e}}ga (private communication).

\begin{prop}\label{prop mult}
Let $k\geq2$ and $\lambda,\mu^0,...,\mu^k\vdash n\geq1$.
We have
$$c^\lambda_{\mu^0,...,\mu^k}(b)=\sum_{\nu\vdash n}c^\lambda_{\mu^0,...\mu^{k-2},\nu}(b)c^\nu_{\mu^{k-1},\mu^k}(b).$$
\begin{proof}
We consider the two functions $\tau^{(k-1)}_b(t,\mathbf{p},\mathbf{q}^{(0)},...,\mathbf{q}^{(k-2)},\mathbf{r})$ and  $\tau^{(1)}_b(t,\mathbf{q}^{(k-1)},\mathbf{q}^{(k)},\mathbf{r})$, and we take their scalar product with respect to the variable $\mathbf{r}$. Since 
$$\langle J_\theta^{(\alpha)}(\mathbf{r}),J_\nu^{(\alpha)}(\mathbf{r})\rangle_\alpha=\delta_{\theta,\nu}j_\theta^{(\alpha)},$$
this scalar product gives the function $\tau^{(k)}_b(t,\mathbf{p},\mathbf{q}^{(0)},...,\mathbf{q}^{(k-2)},\mathbf{q}^{(k-1)},\mathbf{q}^{(k)})$. On the other hand, the expansion of these functions in  power-sum basis can be written:

$$\tau^{(k-1)}_b(t,\mathbf{p},\mathbf{q}^{(0)},...,\mathbf{q}^{(k-2)},\mathbf{r})=\sum_{n\geq0}t^n\sum_{\lambda,\mu^0,...,\mu^{k-2},\nu\vdash n}\frac{c^\lambda_{\mu^0,...,\mu^{k-2},\nu}}{z_\lambda(1+b)^{\ell(\lambda)}}p_\lambda q^{(0)}_{\mu^0}... q^{(k-2)}_{\mu^{k-2}}r_\nu,$$
and 
$$\tau^{(1)}_b(t,\mathbf{q}^{(k-1)},\mathbf{q}^{(k)},\mathbf{r})=\sum_{n\geq0}t^n\sum_{\nu,\mu^{k-1},\mu^k\vdash n}\frac{c^\nu_{\mu^{k-1},\mu^k}}{z_\nu(1+b)^{\ell(\nu)}}  q^{(k-1)}_{\mu^{k-1}}q^{(k)}_{\mu^{k}}r_\nu.$$
We conclude by taking the scalar product of the two last equations.
\end{proof}
\end{prop}
The previous property can be used to extend some results known for coefficients $c$ with three parameters (the case $k=1$) to the general case. In particular, we can deduce the following corollary.

\begin{cor}\label{cor mult}
\Cref{Pos conj} for the coefficients $c^\lambda_{\mu^0,...,\mu^k}$ when $k=1$ implies the conjecture for $c^\lambda_{\mu^0,...,\mu^k}$ for any $k\geq1$.
\begin{proof}
Use induction on $k$ and \cref{prop mult}.
\end{proof}
\end{cor}

\noindent As mentioned in the introduction, the polynomiality of the quantities $c^\lambda_{\mu^0,...,\mu^k}$  has been proved in \cite{DF16} when $k=1$.  This can be generalized for any $k\geq1$. 

\begin{thm}\label{polyc}
For all $\lambda,\mu^0,...,\mu^k\vdash n$, the coefficient $c^\lambda_{\mu^0,..,\mu^k}(b)$ is a polynomial with rational coefficients, and we have the following bounds on the degree:
$$\deg(c^\lambda_{\mu^0,..,\mu^k})\leq \min_{-1\leq i\leq k} d_{i}(\lambda,\mu^0,...,\mu^k),$$
where $$d_{-1}(\lambda,\mu^0,...,\mu^k):=kn+\ell(\lambda)-\big(\ell(\mu^0)+...+\ell(\mu^k)\big),$$
and
$$d_i(\lambda,\mu^0,...,\mu^k):=kn-\sum _{j\neq i }\ell(\mu^j),$$ for  $0\leq i\leq k$.\begin{proof}
The polynomiality and the bound $d_{-1}$ follow from \cite[Proposition B.2]{DF16} and \cref{prop mult}. To deduce the other bounds, we use  the symmetry of coefficients $c^\lambda_{\mu^0...\mu^k}$ in partitions $\mu^i$ for $0\leq i\leq k$ and the following relation that exchanges $\lambda$ and  $\mu^0$ (see \cref{defc}):
\begin{equation}\label{sym2}
    \frac{c^\lambda_{\mu^0...\mu^k}}{z_\lambda(1+b)^{\ell(\lambda)}}=\frac{c^{\mu^0}_{\lambda,\mu^1...\mu^k}}{z_{\mu^0}
    (1+b)^{\ell(\mu^0)}}.
    \qedhere
\end{equation}
\end{proof}
\end{thm}
\noindent In \cref{subsec upperbound}, we give a combinatorial interpretation of the term associated to each one of these bounds. We now state some results that will be useful in \cref{subsec upperbound}.

Using \cite[Lemma 5.7 and Lemma 5.14]{La09}, the polynomiality of the coefficients $c^\lambda_{\mu^0...\mu^k}$, together with \cref{prop mult} we deduce that: 
$$c^\lambda_{\mu^0...\mu^k}=\sum_{0\leq i\leq \left\lfloor\frac{d_{-1}}{2}\right\rfloor}a_i b^{d_{-1}-2i}(1+b)^i,$$ where $a_i\in\mathbb Q$ and $d_{-1}:=d_{-1}(\lambda,\mu^0,..,\mu^k)$. 
The previous equation has the following implication:
\begin{cor}\label{cor b=-1}
For all $\lambda,\mu^0,...,\mu^k\vdash n\geq1$
$$[b^{d_{-1}}]c^\lambda_{\mu^0...\mu^k}=(-1)^{d_{-1}}c^\lambda_{\mu^0...\mu^k}(-1).$$
\end{cor}
\noindent The polynomiality of coefficients $h^\lambda_{\mu^0,...,\mu^k}$ has been deduced from the polynomiality of $c^\lambda_{\mu^0,...,\mu^k}$ when $k=1$ in \cite{DF16}. The proof works in a similar way for $k\geq1$. We obtain the following theorem.
\begin{thm}\label{polyh}
For all $\lambda,\mu^0,...,\mu^k\vdash n\geq1$, the coefficient $h^\lambda_{\mu^0,..,\mu^k}(b)$ is a polynomial in $b$ with rational coefficients, and we have the following bound on its degree:
$$\deg(h^\lambda_{\mu^0,..,\mu^k})\leq  kn+2-\big(\ell(\lambda)+\ell(\mu^0)+...+\ell(\mu^k)\big).$$
\end{thm}
\noindent Using specializations in \cref{defh}, we obtain the following property for the coefficients  $h^\lambda_{\mu^0...\mu^{k}}$ (see \cite[Proposition 4.1]{D17} for the proof in the case $k=1$);
\begin{lem}\label{Somh}
For $\lambda,\mu^0,...,\mu^{k-1}\vdash n$ we have
$$\sum_{\tau\vdash n}h^\lambda_{\mu^0...\mu^{k-1},\tau}(b)=(1+b)^{kn+1-(\ell(\lambda)+\ell(\mu^0)+...+\ell(\mu^{k-1}))}\sum_{\tau\vdash n}h^\lambda_{\mu^0...\mu^{k-1},\tau}(0).$$
\end{lem}
Hence, we deduce the following corollary that will be useful in the proof of \cref{top degree 1}.
\begin{cor}\label{corh}
For $\lambda,\mu^0,...,\mu^{k-1}\vdash n$ we have
$$[b^{kn+1-(\ell(\lambda)+\ell(\mu^0)+...+\ell(\mu^{k-1}))}]h^\lambda_{\mu^0...\mu^k}=\delta_{\mu^k,[n]}\sum_{\tau\vdash n}h^\lambda_{\mu^0...,\mu^{k-1},\tau}(0),$$
where $\delta$ is the Kronecker delta.
\begin{proof}
If $\mu^k\neq[n]$ then from \cref{polyh} we have $$[b^{kn+1-(\ell(\lambda)+\ell(\mu^0)+...+\ell(\mu^{k-1}))}]h^\lambda_{\mu^0...,\mu^{k-1},\mu^k}(b)=0.$$
The previous lemma finishes the proof.
\end{proof}
\end{cor}

\subsection{Top degree in coefficients \texorpdfstring{$c^\lambda_{\mu^0,...,\mu^k}$}{}}\label{subsec upperbound}
\cref{polyc} gives $k+2$ upper bounds on the degrees of coefficients $c^\lambda_{\mu^0,...,\mu^k}$. Using the symmetry property, we can see that the bounds $d_i(\lambda,\mu^0,...,\mu^k)$ for $0\leq i\leq k$ are equivalent.
\cref{top degree 1} gives a combinatorial interpretation for the coefficient in $c^\lambda_{\mu^0,...,\mu^k}$ associated to the bounds  $d_{i}(\lambda,\mu^0,...,\mu^k)$ for $0\leq i\leq k$ and \cref{top degree 2} gives an interpretation for the bound $d_{-1}(\lambda,\mu^0,...,\mu^k)$. The bound $d_{-1}(\lambda,\mu^0,...,\mu^k)$ was investigated in \cite{B21} when $k=1$ and the combinatorial interpretation was given in terms of unhandled maps, while we give here an interpretation with orientable maps with a different proof. In fact, there exists a bijection between the two objects, showing that \Cref{top degree 2} for $k=1$ is equivalent to the result of \cite{B21} (see \cite[Theorem 1.8]{CJS17}).

As explained in \cref{sec gs of const}, the labelling of constellations is simpler in the orientable case. We introduce the following definition of labelling for orientable constellations that will be used to state the main results of this section.

\begin{defi}\label{def or face labelled}
If $\mathbf{M}$ is an orientable $k$-constellation of size $n$:
We say that $\mathbf{M}$:
\begin{itemize}
\item is labelled if its hyperedges are labelled by $\{1,...,n\}$, when we consider $\mathbf{M}$ as a hypermap (see \cref{sec cons}). In terms of right-paths, this is equivalent to label the right-paths of the constellation $\mathbf{M}$ traversed from the corner of color 0 to the corner of color $k$, when $\mathbf{M}$ is equipped with the canonical orientation.
    \item has \textit{labelled faces} if each face has a distinguished corner, and the faces of same size are labelled. 
\end{itemize}

\end{defi}

Note that the definition of face-labelling that we give here for orientable constellations is slightly different from the definition given in \cref{sec face-labelled}; in each face we do not choose an orientation for the distinguished corner. The reason is that in the orientable case all faces have a canonical orientation (see \cref{def const}). We also introduce the following definition.

\begin{defi}\label{def or labelled cc}
Let $\lambda$ be a partition. We say that a $k$-constellation $\mathbf{M}$ is \textit{$\lambda$-connected}, if $\lambda$ is the partition obtained by reordering the sizes of the connected components of $\mathbf{M}$.
We say that a $k$-constellation $\mathbf{M}$ has \textit{labelled connected components}, if each connected component is rooted, and the connected components of the same size are labelled, \textit{i.e.} for $r\geq1$ if $\mathbf{M}$ has $j$ connected components of size $r$, they are labelled with $\{1,..,j\}$.
For every partitions $\lambda,\nu,\mu^0,...,\mu^k\vdash n\geq 1$, we denote $\tilde{h}^{\lambda,\nu}_{\mu^0,...,\mu^k}$ the number of labelled orientable $\nu$-connected $k$-constellations with profile $(\lambda,\mu^0,...,\mu^k)$.
Finally, we say that a $k$-constellation has \textit{partial profile} $(\lambda,\mu^0,...,\mu^{k-1},\bullet)$ if its profile is given by $(\lambda,\mu^0,...,\mu^{k-1},\mu^k)$ for some partition $\mu^k$.
\end{defi}

\begin{thm}\label{top degree 1}
For all $\lambda,\mu^0,...,\mu^k\vdash n\geq1$, the top degree $[b^{d_k}]c^\lambda_{\mu^0\mu^1...\mu^k}$ is equal to the number of $\mu^k$-connected orientable $k$-constellations with labelled faces with partial profile $(\lambda,\mu^0,...\mu^{k-1},\bullet)$, and where $d_k:=d_k(\lambda,\mu^0,...,\mu^k)$.

\begin{proof}
From Equations \eqref{defc}, \eqref{defh} and by developing the exponential in \cref{eqPsi}, we obtain
$$\frac{c^\lambda_{\mu^0... \mu^k}}{z_\lambda(1+b)^{\ell(\lambda)}}=\sum_{r\geq1}\frac{1}{r!}\sum_{(n_i)}\sum_{(\lambda_{(i)},\mu^0_{(i)},...,\mu^k_{(i)})}\prod_{1\leq i\leq r}\frac{h^{\lambda_{(i)}}_{\mu^0_{(i)}... \mu^k_{(i)}}(b)}{n_i(1+b)}
,$$
where the second sum is taken over $r$-tuples of positive integers which sum to $r$, and the third sum is taken over $r$-tuples $(\lambda_{(i)},\mu^1_{(i)},...,\mu^k_{(i)})_{1\leq i\leq r}$, such that $\bigcup\limits_{1\leq i\leq r}\lambda_{(i)}=\lambda$, $\bigcup\limits_{1\leq i\leq r}\mu_{(i)}^j=\mu^j$, for all $j\in\llbracket 0,k\rrbracket$  and $n_i=|\lambda_{(i)}|=|\mu_{(i)}^j|$, for all $i\in\llbracket 1,r\rrbracket$ and $j\in\llbracket 0,k\rrbracket$.\\
The last equality can be rewritten as follows
$$\frac{c^\lambda_{\mu^0... \mu^k}}{z_\lambda}=\sum_{r\geq1}\frac{1}{r!}\sum_{(n_i)}\sum_{(\lambda_{(i)},\mu^0_{(i)},...,\mu^k_{(i)})}\prod_{1\leq i\leq r}\frac{(1+b)^{\ell(\lambda_{(i)})-1}h^{\lambda_{(i)}}_{\mu^0_{(i)}... \mu^k_{(i)}}(b)}{n_i}.$$
But from \cref{polyh} we know that for all $i$, $h^{\lambda_{(i)}}_{\mu^0_{(i)}... \mu^k_{(i)}}(b)$ is a polynomial in $b$ and
$$\deg\left(\frac{(1+b)^{\ell(\lambda_{(i)})-1}h^{\lambda_{(i)}}_{\mu^0_{(i)}... \mu^k_{(i)}}(b)}{n_i}\right)\leq kn_i+1-\sum_{0\leq j\leq k}\ell(\mu^j_{(i)})\leq kn_i-\sum_{0\leq j\leq k-1}\ell(\mu^j_{(i)}).$$
Taking the product over $i$, it gives us  
$$\deg\left(\frac{1}{r!}\prod_{1\leq i\leq r}\frac{(1+b)^{\ell(\lambda_{(i)})-1}h^{\lambda_{(i)}}_{\mu^0_{(i)}... \mu^k_{(i)}}(b)}{n_i}\right)\leq kn-\sum_{0\leq j\leq k-1}\ell(\mu^j)=d_k.$$
To have equality in the last line, we should have  $\ell(\mu^k_{(i)})=1$ for all $i\in \llbracket r \rrbracket$, in other terms $\mu^k_{(i)}=[n_i]$  (and hence $r$ should be equal to $\ell(\mu^k)$). 
Therefore, one has
\begin{align}\label{eqthm1}
    [b^{d_k}]\frac{c^\lambda_{\mu^0,...,\mu^k}}{z_\lambda}
    &=\frac{1}{\ell(\mu^k)!}\sum_{(\lambda_{(i)},\mu^0_{(i)}, ...,n_i)}[b^{d_k}]\prod_{1\leq i\leq \ell(\mu^k)}\frac{(1+b)^{\ell(\lambda_{(i)})-1}h^{\lambda_{(i)}}_{\mu^0_{(i)},...,\mu^{k-1}_{(i)}, [n_i]}}{n_i}\hspace{1cm}\\
    &=\frac{1}{\ell(\mu^k)!}\sum_{(\lambda_{(i)},\mu^0_{(i)}, ...,n_i)}\prod_{1\leq i\leq \ell(\mu^k)}\left[b^{kn_i-\sum\limits_{0\leq j\leq k-1}\ell(\mu^j_{(i)})}\right]\frac{(1+b)^{\ell(\lambda_{(i)})-1}h^{\lambda_{(i)}}_{\mu^0_{(i)},...,\mu^{k-1}_{(i)}, [n_i]}(b)}{n_i},\nonumber
    \end{align}
 where the sums run over $\ell(\mu^{k})$-tuples  $(\lambda_{(i)},\mu^0_{(i)},...,\mu^{k-1}_{(i)},n_i)_{1\leq i\leq \ell(\mu^{k})}$, such that $(n_i)_{1\leq i\leq \ell(\mu^{k})}$ is a reordering of $\mu^{k}$, and $\bigcup\limits_{1\leq i\leq r}\lambda_{(i)}=\lambda$, $\bigcup\limits_{1\leq i\leq r}\mu_{(i)}^j=\mu^j$, for all $j\in\llbracket 0,k-1\rrbracket$. 
 From \cref{corh}, we know that
\begin{equation}\label{eqsomh}
    \left[b^{kn_i-\sum\limits_{0\leq j\leq k-1}\ell(\mu^j_{(i)})}\right](1+b)^{\ell(\lambda_{(i)})-1}h^{\lambda_{(i)}}_{\mu^0_{(i)}... \mu^{k-1}_{(i)}, [n_i]}(b)=\sum_{\tau_{(i)}\vdash n_i} h^{\lambda_{(i)}}_{\mu^0_{(i)}... \mu^{k-1}_{(i)},\tau_{(i)}}(0).
\end{equation}
 Hence $$[b^{d_k}]\frac{c^\lambda_{\mu^0,...,\mu^k}}{z_\lambda}=\frac{1}{\ell(\mu^k)!}\sum_{(\lambda_{(i)},\mu^0_{(i)}, ...\mu^{k-1}_{(i)},\tau_{(i)})}\prod_{1\leq i\leq \ell(\mu^k)}\frac{h^{\lambda_{(i)}}_{\mu^0_{(i)},...\mu^{k-1}_{(i)},\tau_{(i)}}(0)}{n_i}.$$
On the other hand, from \cref{Thm b=0} we know that  $h^{\lambda_{(i)}}_{\mu^0_{(i)},...,\tau_{(i)}}(0)$ is the number of rooted connected orientable $k$-constellations.
Then for every partitions $\lambda_{(i)},\mu^0_{(i)},...,\mu^k_{(i)}\vdash n_i$
$$\frac{h^{\lambda_{(i)}}_{\mu^0_{(i)},...\mu^{k-1}_{(i)},\tau_{(i)}}(0)}{n_i}=\frac{\tilde{h}^{\lambda_{(i)},[n_i]}_{\mu_{(i)}^0,...\mu^{k-1}_{(i)},\tau_{(i)}}}{n_i!},$$
where $\tilde{h}^{\lambda_{(i)},[n_i]}_{\mu_{(i)}^0,...\mu^{k-1}_{(i)},\tau_{(i)}}$ is the quantity defined in \cref{def or labelled cc}. Hence, \cref{eqthm1} can be rewritten as follows
\begin{align*}
     \frac{[b^{d_k}]c^\lambda_{\mu^0,...,\mu^k}}{z_\lambda}=&\frac{\tilde h^{\lambda,\mu^k}_{\mu^0,...\mu^{k-1},\bullet}}{n!}.
\end{align*}
Finally, we multiply $\tilde h^{\lambda,\mu^k}_{\mu^0,...\mu^{k-1},\bullet}$ by $\frac{z_\lambda}{n!}$ to pass from labelled constellations to face-labelled constellations.
\end{proof}
\end{thm}

We now deduce from \cref{top degree 1} an analog theorem for the bound $d_{-1}(\lambda,\mu^0,..,\mu^k)$. 

\begin{thm}\label{top degree 2}
For $\lambda,\mu^0,...,\mu^k\vdash n\geq1$, the top degree term $[b^{d_{-1}}]c^\lambda_{\mu^0,...,\mu^k}$ is equal to the number of $\lambda$-connected orientable $k$-constellation with labelled connected components and partial profile $(\mu^k,\mu^0,...\mu^{k-1},\bullet)$, where $d_{-1}:=d_{-1}(\lambda,\mu^0,..,\mu^k)$. 
\begin{proof}
From $\cref{sym2}$, we know that 
$$[b^{d_{-1}}]c^\lambda_{\mu^0,...,\mu^k}
=\frac{z_\lambda}{z_{\mu^k}}[b^{d_k(\mu^k,\mu^0,...,\mu^{k-1},\lambda)}]c^{\mu^k}_{\mu^0,...,\lambda}.$$
We apply \cref{top degree 1} and  we multiply by $z_\lambda$ to choose the labels of the connected components and we divide by $z_{\mu^k}$ to forget the labels of the faces, which concludes the proof.
\end{proof}
\end{thm}

\noindent {\bf Acknowledgements.} 
The author wishes to thank Guillaume Chapuy and Valentin Féray for suggesting the problem and for many useful discussions.

\bibliographystyle{alpha}
\bibliography{biblio.bib}
\end{document}